\documentclass[11pt, oneside]{article}   	% use "amsart" instead of "article" for AMSLaTeX format
\usepackage{geometry}                		% See geometry.pdf to learn the layout options. There are lots.

\usepackage{graphicx}				% Use pdf, png, jpg, or epsÂ§ with pdflatex; use eps in DVI mode

\usepackage{setspace}
\usepackage[utf8]{inputenc}
\usepackage{amsmath}
\usepackage{mathtools}
\allowdisplaybreaks
\usepackage{amssymb}
\usepackage{todonotes}
\usepackage{xcolor}
\usepackage{longtable}
\usepackage{booktabs}

\usepackage{arydshln}
\makeatletter
\def\adl@drawiv#1#2#3{%
	\hskip.5\tabcolsep
	\xleaders#3{#2.5\@tempdimb #1{1}#2.5\@tempdimb}%
	#2\z@ plus1fil minus1fil\relax
	\hskip.5\tabcolsep}
\newcommand{\cdashlinelr}[1]{%
	\noalign{\vskip\aboverulesep
		\global\let\@dashdrawstore\adl@draw
		\global\let\adl@draw\adl@drawiv}
	\cdashline{#1}
	\noalign{\global\let\adl@draw\@dashdrawstore
		\vskip\belowrulesep}}
\makeatother

\usepackage{siunitx}
\usepackage{pdflscape}
\usepackage{multirow}
\usepackage{colortbl}
\usepackage{natbib}
\usepackage{paralist}
\usepackage[hidelinks]{hyperref}
\usepackage{algpseudocode}
\usepackage{makecell}
\usepackage{subcaption}
\usepackage{threeparttable}
\usepackage{tabularx}
\usepackage{wasysym}
\usepackage{tikz}
\usetikzlibrary{arrows,automata,shapes,calc,positioning}
\usetikzlibrary{decorations.markings}
\usetikzlibrary{decorations.pathreplacing}
\usepackage{forloop} 
\newcounter{ct}

\usepackage{pdflscape}
\usepackage{bbm}

\newcounter{Theo}%[section]
\newcounter{Algo}%[section]
\newcounter{Assu}%[section]

\newtheorem{Theorem}[Theo] {Theorem}

\newtheorem{Example}[Theo] {Example}

\newcommand{\id}{1\hspace{-0,9ex}1}
\usepackage{natbib}
 \bibpunct[, ]{(}{)}{,}{a}{}{,}%
 \def\bibfont{\small}%
\usepackage{geometry}
\title{A Tight Formulation for the Dial-a-Ride Problem}
\usepackage{authblk}
\author[1]{Daniela Gaul }
\author[1]{Kathrin Klamroth}
\author[2]{Christian Pfeiffer}
\author[2]{Arne Schulz}
\author[1]{Michael Stiglmayr}
\affil[1]{School of Mathematics and Natural Sciences, University of Wuppertal, Gau{\ss}stra{\ss}e 20, Wuppertal 42119, Germany, \textit{\{gaul,klamroth,stiglmayr\}@math.uni-wuppertal.de}}
\affil[2]{Institute of Operations Management, Universit\"at Hamburg, Moorweidenstra{\ss}e 18, 20148 Hamburg, Germany, \textit{\{christian.pfeiffer,arne.schulz\}@uni-hamburg.de} \newline}

\begin{document}
	
	\newpage\null\thispagestyle{empty}
	\begin{center}
		\textbf{\Large Remark}\\
	\end{center}

	Theorem 2 of this paper is false. Therefore it has been removed in the published version to be found at https://doi.org/10.1016/j.ejor.2024.09.028.
	\newpage
	
	\maketitle
	
	\begin{abstract}
Ridepooling services play an increasingly important role in modern transportation systems. With soaring demand and growing fleet sizes, the underlying route planning problems become increasingly challenging. In this context, we consider the dial-a-ride problem (DARP): Given a set of transportation requests with pick-up and delivery locations, passenger numbers, time windows, and maximum ride times, an optimal routing for a fleet of vehicles, including an optimized passenger assignment, needs to be determined. 

We present tight mixed-integer linear programming (MILP) formulations for the DARP by combining two state-of-the-art models into novel \emph{location-augmented-event-based} formulations. Strong valid inequalities and lower and upper bounding techniques are derived to further improve the formulations. We then demonstrate the theoretical and computational superiority of the new model: First, the formulation is tight in the sense that, if time windows shrink to a single point in time, the linear programming relaxation yields integer (and hence optimal) solutions. Second, extensive numerical experiments on benchmark instances show that computational times are on average reduced by 49.7\% compared to state-of-the-art event-based approaches.
	\end{abstract}


\maketitle
%%%%%%%%%%%%%%%%%%%%%%%%%%%%%%%%%%%%%%%%%%%%%%%%%%%%%%%%%%%%%%%%%%%%%%

\section{Introduction}
The dial-a-ride problem (DARP) is a well-studied routing problem \citep[e.g.,][]{HSKLPT18}. It has topical applications in modern transportation systems. While ridepooling  services  are probably the most important area of application \citep[e.g.,][]{Fol20,PS22}, DARPs are also relevant for the transportation of patients \citep[e.g.,][]{MBCB17} or, in a related branch of research, for the transportation of goods. %in the related pick-up-and-delivery problem with time windows (PDPTW). 
Given a fixed number of vehicles, the goal is to identify optimal vehicle routes together with an optimized passenger assignment so that a set of customer requests is served in a best possible way. A customer request typically consists of a number of passengers that want to be transported from a pick-up to a delivery location within specified time windows and within a maximum acceptable ride time. While several requests may share the same vehicle for parts of their trips, the capacity of a vehicle may not be exceeded at any point in time.

This paper is motivated by the ridepooling service HolMichApp (\url{https://www.holmich-app.de}) which was launched in 2019 in the mid-sized city of Wuppertal in Germany with the purpose of improving mobility and reducing congestion and the use of private cars.  %Like HolMichApp, a large number of providers have launched on-demand ridepooling services for the general public, especially in recent years \citep{Fol20}.
Ridepooling services play an important role in the transformation to a sustainable future by improving mobility and reducing congestion and emissions. %Because of that, e.g., the German government funds several ridepooling services as part of a sustainable public transport system \citep{BMDV22}. Beside public providers which work as part of the \ms{local public transport service}, we also find private providers in the ridepooling market (like MOIA in Hamburg (\url{https://www.moia.io/en})). 
To make ridepooling services competitive and efficient, optimized routing and customer assignments are of central importance. Hence, the exact and fast solution of realistically-sized DARP instances in a rolling horizon environment is pivotal for the success of ridepooling services in particular, and for a sustainable improvement of public transportation systems in general.

In this paper, we present new mixed-integer linear programming (MILP) formulations for the static DARP that build upon state-of-the-art event-based (EB) and location-based (LB) models, using the best of both worlds. Our results transfer to the dynamic DARP in two ways: First, \citet{GKS21} showed that EB formulations can be applied efficiently in a rolling-horizon framework for the dynamic DARP. This also holds for the new models developed in this paper. Second, we show in this paper that the new models are tighter than well-known LB models (see Theorem \ref{T1}) and yield an integral polyhedron if the time windows fulfill additional conditions (see Theorem \ref{Net}). This implies that integer optimal solutions can be determined as optimal solutions of the linear programming (LP) relaxation within milliseconds.
%leads to integer optimal solutions of the linear programming (LP) relaxation if time windows reduce to single points in time. 
We note that \citet{BJM19} used a similar result for the large-scale dynamic optimization of taxis. %We describe the contribution of the paper in detail after discussing the relevant literature.

\subsection{Related Work}
Since this paper considers mixed-integer linear programming formulations for the DARP, %and their optimal solution, 
the following short literature review is focused on exact methods for DARP and for the closely related pick-up and delivery problem with time windows (PDPTW). For a broader introduction  we refer to the surveys of \citet{MBC17} and \citet{HSKLPT18}. 

Most exact methods are based on branch-and-cut or branch-and-price frameworks that are tailored for different MILP formulations. \citet{Cor06} presented a 3-index MILP formulation for the DARP where binary variables indicate whether a certain vehicle $k$ drives directly from location $i$ to location $j$. %\ms{Moreover,} \citet{Cor06} introduced several classes of valid inequalities and a set of benchmark instances. 
One year later, \citet{RCL07} reformulated the model of \citet{Cor06} by omitting the index for the vehicles. Now, a binary variable indicates whether any vehicle drives directly from location $i$ to location $j$. Both models are typical examples of \emph{location-based (LB) formulations} in which binary variables indicate whether a vehicle drives directly between two locations. \citet{Cor06} and \citet{RCL07} suggest a branch-and-cut framework for the LB models, and \citet{RCL07} derived additional valid inequalities. Particularly the fork constraints and the reachability constraints that were adapted from \citet{L06} proved to be very effective. LB models are discussed in more detail in Section~\ref{sec:lbMILP}. %we present one of the two models from \citet{RCL07} \ms{in more detail}.

In branch-and-price methods, the DARP is usually subdivided into two subproblems: The master problem where tours are selected, and the pricing problem where the tours are generated. \citet{RC09} used a branch-and-cut-and-price approach for the PDPTW. %pick-up and delivery problem with time windows. 
They compared two formulations for the pricing problem, an elementary and a non-elementary shortest path problem. 
The branch-and-cut-and-price algorithm suggested in \citet{GI15} instead uses a pricing problem that includes ride time constraints explicitly, making it more difficult to solve. 

Comparing branch-and-cut and branch-and-price, we observe that the LB formulation has a quadratic number of binary variables  (in the number of locations). In the worst case, all feasible tours need to be constructed in the branch-and-price approach. However, the LB formulation requires a potentially exponential number of constraints to ensure that corresponding pick-up and delivery nodes are visited in the same tour. This is not the case in the pricing problem, as it can be formulated as a resource constrained shortest path problem. In conclusion, different formulations have a trade-off between the number of binary variables and the number of required constraints.

In recent years, new approaches found different compromises with respect to this trade-off. %\kk{i.e., compromises between} %which can be classified in between of the two extremes given by 
%the LB formulation and the branch-and-price approach. 
\citet{RF21} introduced route fragments which are segments of a vehicle tour such that the vehicle is empty at the beginning and at the end of a fragment. Thus, vehicle tours can be compounded by fragments. \citet{RF21} used the fact that time windows typically strongly reduce the number of possible fragments. In their approach, binary decision variables indicate whether a fragment is used in the solution or not and whether fragments are scheduled in direct succession within a tour. They used a branch-and-cut framework to solve their model. \citet{GKS22} introduced an \emph{event-based} (EB) formulation where an event represents the current occupancy of the vehicle. In this formulation, binary variables indicate whether two events occur in direct succession in a vehicle's tour. \citet{GKS22} used the fact that the vehicle capacity is typically rather small in practice, as the number of events strongly depends on the vehicle capacity. We review their approach in detail in Section~\ref{sec:ebMILP}. Both aforementioned approaches by \citet{RF21} and \citet{GKS22} have the advantage that the vehicle capacity constraints are enforced by the choice of fragments and events, respectively. Moreover, precedence constraints are either ensured by the choice of fragments or can be implemented implicitly by the selection of arcs connecting events. This considerably reduces the number of constraints in EB models as compared to LB formulations.

\subsection{Our Contribution} \label{Contribution}
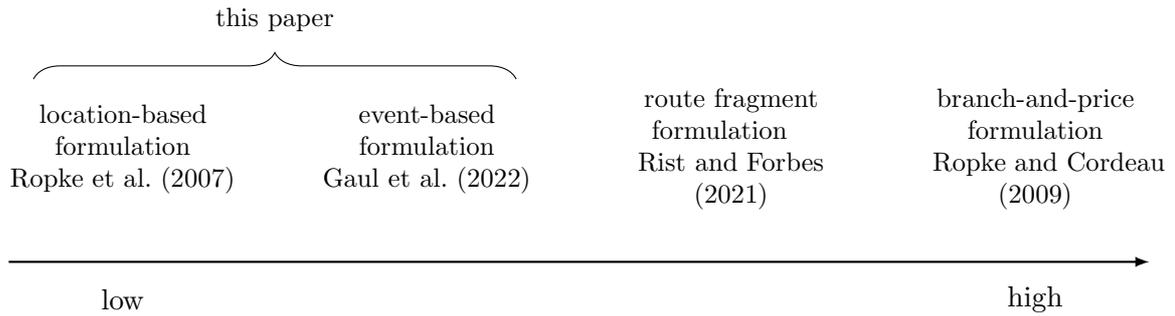
\begin{figure}[htbp]
	\centering
	\begin{tikzpicture}
\draw[thick, -latex] (-1.5,0) -- (13.5,0);
\node at (0,-0.5) {low};
\node at (12,-0.5) {high};
%\node at (6,-0.5) {information given by a binary variable};
\node (A) at (0,1.5) {\begin{minipage}{8em}\small\centering location-based formulation \citet{RCL07} \end{minipage}};
\node (B) at (4,1.5) {\begin{minipage}{8em} \small\centering event-based formulation \citet{GKS22} \end{minipage}};
\node at (8,1.5) {\begin{minipage}{7em} \small\centering route fragment formulation \newline\citet{RF21} \end{minipage}};
\node at (12,1.5) {\begin{minipage}{7em} \small\centering branch-and-price formulation \citet{RC09} \end{minipage}};

\draw[decorate,decoration={brace,raise=2mm,amplitude=10pt}] (A.north west)+(0.5,0) to node[above=7mm] {\small this paper} ([shift=({-0.5,0})]B.north east);

\end{tikzpicture}
	\caption{Classification of formulations according to the information given by a binary variable with value~1}
	\label{fig:class}
\end{figure}

The above literature review implies a classification of MILP formulations with regard to the information contained in the binary variables, see Figure~\ref{fig:class} for an illustration. In this paper, we propose novel MILP formulations that combine the advantages from the LB formulation of \citet{RCL07} and the EB model of \citet{GKS22}. The new models are further improved by adapting problem specific valid inequalities inspired by \citet{SP22}.
%in the nodes of the branch-and-bound tree to improve the search in the LB formulation. 
Our main contributions are:
\begin{itemize}
    \item We introduce a \emph{location-augmented-event-based (LAEB) formulation} and an \emph{aggregated location-augmented-event-based (ALAEB) formulation}. The ALAEB is a slightly adapted variant of the LAEB that has significantly fewer binary variables and thus potentially reduces the size of the branch-and-bound tree.
		\item We proof that both new formulations have a provably tighter LP relaxation than previous LB models (see Theorem \ref{T1}).
    \item We show that both new formulations are tight in the sense that their LP relaxation have integral solutions if the time windows fulfill additional conditions (see Theorem \ref{Net}). This tightness is also reflected by an average root node gap of only 1.6\% in our computational study (for general time windows).
	\item Preprocessing the event-based graph reduces the size of both new MILP formulations by 22\% on average.	
	\item New valid inequalities are derived for the EB, LAEB, and ALAEB formulations. 
	\item In a comprehensive computational study,	the LAEB formulation leads to a 49.7\% reduction of computational time on benchmark instances as compared to state-of-the-art EB formulations.
	\item Our formulation is hence suitable for a rolling-horizon framework and for practical applications of a ridepooling service as HolMichApp.
    %\item \kk{We introduce a \emph{location-augmented-event-based (LAEB) formulation} that has a provably tighter LP relaxation than previous LB models.}
    %\item We show that \kk{the LAEB formulation is tight in the sense that its LP relaxation has integral solutions} if the time windows shrink to a single point in time. 
    %\kk{This tightness is also reflected by an average root node gap of only 1.6\% in our computational study (for general time windows). }
    %\todo[inline]{@Arne: Gilt auch allgemeiner wenn paarweise Reihenfolge festgelegt ist?}
    %\item \kk{An \emph{aggregated location-augmented-event-based (ALAEB) formulation} is presented as an alternative model that has significantly fewer binary variables and thus potentially reduces the size of the branch-and-bound tree.} 
    %\item \kk{Preprocessing the event-based graph reduces the size of both new MILP formulations  by 22\% on the average.} 	
	%\item New valid inequalities are derived for the EB, LAEB\kk{, and ALAEB} formulations. 
	%\item \kk{In a comprehensive computational study,	the LAEB formulation leads to a 49.7\% reduction of computational time on benchmark instances as compared to state-of-the-art EB formulations.}
	%\item Our formulation is hence suitable for a rolling-horizon framework and for practical applications of a ridepooling service as HolMichApp.
\end{itemize}
The paper is structured as follows: First, we recapitulate the LB formulation by \citet{RCL07} and the EB formulation by \citet{GKS22} in detail and present our new formulations in Section~\ref{sec:MILPs}. Our theoretical investigation of the new formulations is presented in Section~\ref{sec:TA}. In Section~\ref{sec:BC}, we introduce methods to improve the performance of branch-and-bound search. The new formulations as well as the introduced methods are evaluated in the computational study in Section~\ref{sec:CS}. Finally, the paper closes with a conclusion in Section~\ref{sec:Con}.

\section{MILP formulations} \label{sec:MILPs}
We consider a DARP with $n$ customers. Each customer request $i \in R = \{1, \ldots, n\}$ contains a number of persons $q_i \in \mathbb{N}\setminus\{0\}$ and a service time $s_i \geq 0$, which can be interpreted as the time the customer needs to enter or exit the vehicle. They have to be transported by one of $K$ identical vehicles with capacity $Q$ between their pick-up location $i^+$ and delivery location $i^-$. The set of pick-up locations is denoted by $P = \{1^+, \ldots,n^+\}$ and the set of delivery locations by $D = \{1^-, \ldots,n^-\}$. The vehicle depot is represented by $0$. To simplify the notation, we set $s_{i^+} = s_{i^-} = s_i$, $s_0 = 0$, $q_{i^+} = q_{i^-} = q_i$ as well as $q_0 = 0$. A customer is not allowed to be in the vehicle for a longer time than a specified maximum ride time $L_i$. For the travel time $\bar{t}_{ij}$ and corresponding routing costs $\bar{c}_{ij}$ between each pair of locations $i, j \in J = P \cup D \cup \{0\}$ we assume that $\bar{c}_{ij}$ and $\bar{t}_{ij}$ are non-negative and fulfill the triangle inequality. Let $\bar{t}_i$ denote the travel time for a direct ride from pick-up to delivery location of request $i$, i.e., $\bar{t}_i = \bar{t}_{i^+i^-}$. W.l.o.g. we assume $0<\bar{t}_i + s_i<L_i$ %$\bar{t}_i + s_i >0$ 
for all $i\in R$. All locations $i\in J$ have to be visited within a given time window $[e_i,\ell_i]$. The beginning of service is given by the lower bound $e_0 = 0$ of the time window at the depot, and there is a fixed duration of service $T$, so that the end of service is given by $\ell_0 \coloneqq e_0 + T$. We assume that $e_{i^+} \geq \bar{t}_{0i^+}$, i.e., each pick-up location can be reached at the beginning of the time window when starting at time 0 in the depot. The objective is to minimize the total routing costs. An overview of the used notation can be found in the electronic companion \ref{Notation}. Note, that we use a bar over parameters and variables if they are indexed by locations and the same type of parameter and variable (without bar), respectively, is used for events.

The DARP can be represented by different graph-based formulations. We first present the LB formulation by \citet{RCL07} in the following section before we present the EB formulation by \citet{GKS22}. Afterwards, we combine both formulations into two location-augmented-event-based formulations in Sections~\ref{IEBM} and \ref{SELF}, respectively, by using the advantages of both the event- and the location-based formulation.

\subsection{Location-Based MILP Formulation} \label{sec:lbMILP}
In this section, we present the LB formulation by \citet{RCL07} which will be used as reference LB model in the following (the first one in their paper; additionally the number of vehicles is restricted). In this model, the vehicle depot is denoted by $0^+$ and $0^-$, raising the possibility of a different start and end depot, however, in this paper, $0^+$ and $0^-$ represent the same vehicle depot $0$ and are only used to stick to the notation from the original paper by \cite{RCL07}. Thus, we set $s_{0^+} = s_{0^-} = s_0$, $q_{0^+} = q_{0^-} = q_0$, $e_{0^+} = e_{0^-} = e_0$, $\ell_{0^+} = \ell_{0^-} = e_0$ as well as $\bar{c}_{0^+j} = \bar{c}_{0^-j} = \bar{c}_{0j}$, $\bar{c}_{i0^+} = \bar{c}_{i0^-} = \bar{c}_{i0}$, $\bar{t}_{0^+j} = \bar{t}_{0^-j} = \bar{t}_{0j}$, and $\bar{t}_{i0^+} = \bar{t}_{i0^-} = \bar{t}_{i0}$. In the LB formulation, each node in the underlying graph is associated with a location $j \in \bar{J} = P \cup D \cup \{0^+, 0^-\}$. The graph is complete and directed, i.e., travelling the arc $(i,j)$ means that the vehicle drives from location $i$ to $j$, where $i,j \in \bar{J}$.

A binary variable $\bar{x}_{ij}$ is associated with each arc is 1 if a vehicle drives directly from location $i$ to location $j$, i.e., the arc $(i,j)$ is used, and 0 otherwise. Moreover, $Q_j$ represents the number of customers in the vehicle after leaving location $j$ and $\bar{B}_j$ is the time when service at location $j$ begins.
In the model, $\mathcal{S}$ is the set of all sets $S$ with $0^+ \in S$,  $0^- \notin S$, and there is at least one request $i$ for which the delivery node is in $S$ but not the pick-up node. This means $\mathcal{S} = \{ S \colon 0^+ \in S \wedge 0^- \notin S \wedge \exists i: (i^+ \notin S \wedge i^- \in S) \}$. Moreover, $\bar{M}_{ij}$ with $i,j\in \bar{J}$ is a sufficiently large constant. Then, the model is formulated as:
\begin{subequations}
	\begin{align}
		Z_{LB} = \min\;& \sum_{i,j \in \bar{J}} \bar{c}_{ij} \cdot \bar{x}_{ij} \label{eq:objFunc}\\
		\text{s.t.}\;&\sum_{i \in \bar{J}} \bar{x}_{ij} = 1\quad&&\forall j\in P\cup D\label{eq:degIn}\\ 
		&\sum_{j \in \bar{J}} \bar{x}_{ij} = 1\quad&&\forall i\in P\cup D\label{eq:degOut}\\
		&\sum_{j\in P} \bar{x}_{0^+j}\leq K&&\label{eq:depot}\\
		&\sum_{i,j\in S} \bar{x}_{ij}\leq |S|-2\quad &&\forall S\in \mathcal{S}\label{eq:SEC}\\
		&\bar{B}_i + s_i + \bar{t}_{ij} - \bar{M}_{ij}\, (1 - \bar{x}_{ij}) \leq \bar{B}_j &&\forall i,j \in \bar{J}\label{eq:time}\\
		&Q_i + \mathcal{I}_{\{ j \in P \}}\, q_j - Q\, (1 - \bar{x}_{ij}) \leq Q_j &&\forall i,j \in \bar{J}\label{eq:cap}\\
		&e_j \leq \bar{B}_j\leq \ell_j&&\forall j\in \bar{J}\label{eq:B}\\
		&\max\{0, \mathcal{I}_{\{ j \in P\}}q_j \}\leq Q_j\leq\min\{Q,Q+ \mathcal{I}_{\{ j \in P\}}q_j \}&&\forall j\in \bar{J} \label{eq:q}\\
		&\bar{B}_{i^-} - \bar{B}_{i^+} - s_i \leq L_i &&\forall i \in P \label{eq:mrt} \\
		&\bar{x}_{ij}\in\{0,1\}&&\forall i,j \in \bar{J}\label{eq:xbin}
	\end{align}
\end{subequations}

Objective function \eqref{eq:objFunc} minimizes the total routing costs. Constraints \eqref{eq:degIn} and \eqref{eq:degOut} ensure that each customer location is visited and left exactly once. Due to constraint \eqref{eq:depot} the depot is left at most $K$ times, i.e., up to $K$ vehicles are used. Because of constraints \eqref{eq:SEC} we are sure that the pick-up and the delivery location of each customer are visited in the same tour and in the correct order. Constraints \eqref{eq:time} and \eqref{eq:B} take care that each location is visited within its time window. Constraints \eqref{eq:cap} model the vehicle load which is not allowed to exceed the vehicle's capacity because of constraints \eqref{eq:q}, where 
\begin{equation*}
	\mathcal{I}_{\{j\in P\}} = 
    \begin{cases}
		1 & j\in P, \\
		-1 & \text{otherwise}.
	\end{cases}
\end{equation*}
No customer is allowed to be longer in the vehicle than their maximum ride time $L_i$, which is ensured by constraints \eqref{eq:mrt}. Finally, Constraints \eqref{eq:xbin} are the binary constraints for the sequence variables. 

As can be seen the size of the set $\mathcal{S}$ grows exponentially in the number of customers $n$. Moreover, constraints \eqref{eq:time} and \eqref{eq:cap} are numerically unfavorable big-$M$-constraints (the $Q$ in \eqref{eq:cap} is also a large constant). The larger $M_{ij}$ is, the closer can a fractional value of $\bar{x}_{ij}$ be to 1 without being binding in constraints \eqref{eq:time} \citep{CRT90,CF06}. 

% The EB formulation presented in the following section overcomes two of the aforementioned drawbacks, as it does not contain \kk{an exponentially growing number of} constraints as constraints \eqref{eq:SEC}, since it relies on a graph structure which \kk{enforces} pairing and precedence constraints implicitly. Also vehicle capacity constraints are incorporated implicitly, i.e., the model does not \kk{require} constraints \eqref{eq:cap}.

\subsection{Event-Based MILP Formulation} \label{sec:ebMILP}
The EB formulation introduced by \citet{GKS22} overcomes two of the aforementioned drawbacks of LB, as it does not contain an exponentially growing number of constraints as constraints \eqref{eq:SEC}, since it relies on a graph structure which enforces pairing and precedence constraints implicitly. Also vehicle capacity constraints \eqref{eq:cap} are incorporated implicitly.
% In this section, we present the EB formulation introduced by \citet{GKS22}. 
In the event-based graph, nodes represent states or events, which are not only associated with the current location of the vehicle but also with the entire state of the vehicle, including the information on all customers sitting in the vehicle.

As at most $Q$ customer requests are transported at the same time by a vehicle (if $q_i = 1$ for all of them), the allocation of customers to a vehicle (and hence the nodes of the event-based graph) can be represented by a $Q$-tuple. If less than $Q$ requests are allocated to the vehicle, the remaining places are marked by 0. This holds also if a single request already contains $Q$ persons. Thus, a zero component in $Q$ does not mean that the vehicle has empty seats. Since all customers present in the vehicle are part of the event description, we can easily ensure that capacity constraints \eqref{eq:cap} are fulfilled by deleting all nodes where the sum of the $q_i$ values of all customers present in the vehicle exceeds $Q$. Thus, we can never reach a situation where the capacity constraint is violated. 

Reordering the customers in the vehicle leads to different representations of the same state. To avoid ambiguities,  we assume that the customers are noted in decreasing order of their index and the event is filled up with zeros on the remaining places. An exception is the first position of the $Q$-tuple, which indicates the current position of the vehicle, i.e., the last location visited. By this, there is a direct mapping of events to locations. On the other hand, each location may be associated with different events with different customers present in the vehicle at the time the location is reached. 

As an example, event $(3^+,2,1)$ represents for a vehicle capacity of $Q = 3$ the situation that the vehicle is at the pick-up location of customer 3 ($3^+$) and customers 1 and 2 are sitting in the vehicle. Another example is the event $(1^-,2,0)$. Here the vehicle delivered customer 1 and customer 2 is still sitting in the vehicle. 

\citet{Cor06} pointed out that time window and ride time constraints can be used to identify incompatible pairs of users. Given users $i, j \in R$, the following reductions are possible \citep[cf.][]{GKS22}:
\begin{itemize}
	\item if it is infeasible to visit locations $ i^+, i^-, j^+, j^-$ both in the order $ j^+ \rightarrow i^+\rightarrow j^- \rightarrow i^-$ and $ j^+ \rightarrow i^+\rightarrow i^- \rightarrow j^-$,  the set $\{(v_1, v_2, \ldots, v_Q)\in V\colon v_1 = i^+, j\in \{v_2, \ldots, v_Q\}\}$ can be removed from the node set.
	\item if it is infeasible to visit locations $ i^+, i^-, j^+, j^-$ both in the order $j^+ \rightarrow i^+ \rightarrow i^- \rightarrow j^-$ and $i^+ \rightarrow j^+ \rightarrow i^- \rightarrow j^-$, the set $\{(v_1, v_2, \ldots, v_Q)\in V\colon v_1 = i^-, j\in \{v_2, \ldots, v_Q\}\}$ can be removed from the node set.
\end{itemize}
Now, let $f^1_{ij}$ and $f^2_{ij} \in\{0,1\}$ encode the feasibility of the paths $j^+ \rightarrow i^+ \rightarrow j^- \rightarrow i^-$ and $j^+ \rightarrow i^+ \rightarrow i^- \rightarrow j^-$, respectively, w.r.t.\ ride time and time window constraints. To simplify the notation, let $f^1_{i0} = f^1_{0i} = f^2_{i0} = f^2_{0i} = 1$. By this, the number of $Q$-tuples representing feasible user allocations can be reduced significantly by computing $f^1_{ij}, f^2_{ij}$ for all pairs of users $i,j\in R$.
Then, the set of pick-up nodes of user $i\in R$ can be expressed as
\begin{equation*}
	\begin{split}
		V_{i^+} \coloneqq \Biggl\lbrace  (v_1, v_2, \ldots, v_Q) \colon  v_1 &= i^+,\; v_l \in R \cup \{0\} \setminus \{i\}, \; f^1_{i,v_l} + f^2_{i,v_l} \geq 1\; \forall l \in\{2, \ldots,Q\},\\
		&\Bigl(v_l > v_{l+1} \vee v_{l+1}=0 \Bigr) \; \forall l \in\{2, \ldots,Q-1\}, \; \sum_{l=1}^Q q_{v_l} \leq Q\Biggr\rbrace.
	\end{split}
\end{equation*}
Analogously, the set of delivery nodes corresponding to $i\in R$ can be expressed as
\begin{equation*}
	\begin{split}
		V_{i^-} \coloneqq \Biggl\lbrace  (v_1, v_2, \ldots, v_Q) \colon  v_1 &= i^-,\; v_l \in R \cup \{0\} \setminus \{i\}, \; f^1_{v_l,i} + f^2_{i,v_l} \geq 1\; \forall l \in\{2, \ldots,Q\},\\
		&\Bigl(v_l > v_{l+1} \vee v_{l+1}=0 \Bigr) \; \forall l \in\{2, \ldots,Q-1\}, \; \sum_{l=1}^Q q_{v_l} \leq Q\Biggr\rbrace.
	\end{split}
\end{equation*}

The depot is represented by the single event $(0, \ldots,0)$. We set $V_0 = \{ (0, \ldots,0) \}$. In total, the nodes of the event-based graph are given by
\[
V = V_0 \cup \bigcup_{i=1}^n V_{i^+} \cup \bigcup_{i=1}^n V_{i^-}.
\]
The arcs correspond to feasible transitions between events. In the aforementioned example, there is no arc between events $(3^+,2,1)$ and $(1^-,2,0)$, as customer 3 is present in the vehicle in the first event but not in the second. However,
\[
(3^+,2,1) \rightarrow (3^-,2,1) \rightarrow (1^-,2,0)
\]
may be feasible transitions. In the left event, customer 3 enters the vehicle with customers 1 and 2 sitting in it at customer 3's pick-up location. In the second, the vehicle is at customers 3's delivery location, i.e., customer 3 leaves the vehicle while customers 1 and 2 are still in it. Then, the vehicle drives to customer 1's delivery location while customer 2 is still sitting in the vehicle. In comparison to the LB approach, arcs still correspond to a transition between two locations. However, while an  arc in the location-based graph represents only that the vehicle drives from one location to another, the information on the other customers sitting in the vehicle is also encoded in the event-based representation. In other words, the LB arc between locations $i$ and $j$ aggregates all event-based arcs between nodes $v$ and $w$ representing these two locations, i.e., with $v_1 = i$ and $w_1 = j$ where $v_1$ and $w_1$ represent the first entries of $v$ and $w$, respectively.
%($w_1$) represents the first entry of $v$ ($w$). 

In line with this connection between both approaches, the set of arcs is
\[
A = \bigcup_{i=1}^6 A_i,
\]
where the subsets $A_i$ are defined as follows:
\begin{itemize}
	\item $A_1$ defines the transition between a customer's pick-up location and a customer's delivery location, i.e.
	\[
	A_1 = \left\{ ((i^+,v_2, \ldots,v_Q),(j^-,w_2, \ldots,w_Q)) \in V \times V: \{ i,v_2, \ldots,v_Q \} = \{ j,w_2, \ldots,w_Q \} \right\}.
	\]
	\item $A_2$ represents transitions between a customer's pick-up location and another customer's pick-up location, i.e.
	\[
	A_2 = \left\{ ((i^+,v_2, \ldots,v_{Q-1},0),(j^+,w_2, \ldots,w_Q)) \in V \times V: \{ i,v_2, \ldots, V_{Q-1}\} = \{ w_2, \ldots, w_Q \} \right\}.
	\]
	\item $A_3$ is the set of all transitions between a customer's delivery location and another customer's pick-up location, i.e.
	\[
	A_3 = \left\{ ((i^-,v_2, \ldots,v_Q),(j^+,v_2, \ldots,v_Q)) \in V \times V: i \neq j \right\}.
	\]
	\item $A_4$ defines the transition between the delivery locations of two customers, i.e.
	\[
	A_4 = \left\{ ((i^-,v_2, \ldots,v_Q),(j^-,w_2, \ldots,w_{Q-1},0)) \in V \times V: \{ v_2, \ldots,v_Q\} = \{j,w_2, \ldots,w_{Q-1}\} \right\}.
	\]
	\item $A_5$ is the set of all transitions between a delivery location and the depot, i.e.
	\[
	A_5 = \left\{ ((i^-,0, \ldots,0),(0, \ldots,0)) \in V \times V \right\}.
	\]
	\item $A_6$, finally, describes the transition between the depot and a pick-up location, i.e.
	\[
	A_6 = \left\{ ((0, \ldots,0),(i^+,0, \ldots,0)) \in V \times V \right\}.
	\]
\end{itemize}
\citet{GKS22} show that the number of nodes and arcs of the event-based graph can be bounded for $n \geq Q+1$ by $\mathcal{O} \left( n^Q \right)$ and $\mathcal{O} \left( n^{Q+1} \right)$, respectively. 
% by $\mathcal{O} \left( n^Q \right)$ for $n \geq Q$ and $\mathcal{O} \left( n^{Q+1} \right)$ for $n \geq Q+1$, respectively.
In comparison, the LB approach requires one node per location, i.e., $\mathcal{O}(n)$, and one arc for each pair of locations, i.e., $\mathcal{O}\left( n^2 \right)$. Thus, the advantage of the LB formulation is the smaller size of the underlying graph and thereby of the number of variables. However, in practical applications $Q$ is usually rather small (typically $Q = 3$ or $Q = 6$), and we can reduce the size of the event-based graph due to time windows, travel and ride times (which is presented in Section~\ref{sec:BC}). On the other hand, the event-based graph has the already mentioned advantage that capacity, pairing, and precedence constraints are modelled implicitly in the graph. 

In the following, we formulate the EB model by \citet{GKS22}. To simplify the notation, let $a = (v,w) \in A$ be the arc between the nodes representing the events $v$ and $w$. Then, routing costs and travel times between locations $v_1$ and $w_1$ can be represented as $c_a$ and $t_a$, respectively. Moreover, let $\delta^{in}(v) = \{(w,v) \in A \}$ and $\delta^{out}(v) = \{(v,w) \in A \}$ be the sets of incoming and outgoing arcs for node $v$. Analogous to the LB model, variable $x_a$ is 1 if arc $a$ is used in the solution and $B_v$ reports the beginning of service at location $v_1$. 

\citet{GKS22} present two model formulations. The following one is the second which led to significantly shorter computation times in their computational evaluation. In some applications of the DARP, we can distinguish between inbound requests, i.e.\ requests that submit a time for pick-up, and outbound requests, i.e. requests that submit a time for delivery, see \cite{Cor06}. In the EB formulation, \cite{GKS22} assume that all pick-up time windows of inbound requests and all delivery time-windows of outbound requests have a length of $TW$, i.e.\ $TW = \ell_j - e_j$.
\begin{subequations}
	\begin{align}
		Z_{EB} = \min\;& \sum_{a \in A} c_a \cdot x_a \label{E:objective}\\
		\text{s.t.}\;&\sum_{a \in \delta^{in}(v)} x_a - \sum_{a \in \delta^{out}(v)} x_a = 0\quad&&\forall v \in V\label{E:flow}\\
		&\sum_{a \in \delta^{in}(v), v \in V_{i^+}} x_a = 1\quad&&\forall i \in R\label{E:pick}\\
		&\sum_{a \in \delta^{out}(0, \ldots,0)} x_a \leq K\quad&&\label{E:K}\\
		&B_w \geq B_v + s_{v_1} + t_{(v,w)} - M_{vw} (1 - x_{(v,w)})\quad&&\forall (v,w) \in A, v \neq (0, \ldots,0)\label{E:time}\\
		&B_w \geq e_0 + t_{(v,w)}x_{(v,w)}\quad&&\forall((0, \ldots,0),w) \in A\label{E:time1}\\
		&e_0 \leq B_{(0, \ldots,0)} \leq \ell_0\quad&&\label{E:tw1}\\
		&e_{i^+} + TW \Biggl( 1 - \sum_{a \in \delta^{in}(v)} x_a \Biggr) \leq B_v \leq \ell_{i^+}\quad&&\forall i \in R, v \in V_{i^+}\label{E:tw2}\\
		&e_{i^-} \leq B_v \leq e_{i^+} + L_i + s_{i^+} + TW \sum_{a \in \delta^{in}(v)} x_a\quad&&\forall i \in R, v \in V_{i^-}\label{E:tw3}\\
		&B_w - B_v - s_{i^+} \leq L_i\quad&&\forall i \in R, v \in V_{i^+}, w \in V_{i^-}\label{E:mrt}\\
		&x_a \in \{0,1\}\quad&&\forall a \in A\label{E:binary} %\\
		%&B_v \geq 0\quad&&\forall v \in V\label{E:Bnn}
	\end{align}
\end{subequations}

   The objective function corresponds to that of the LB formulation; in both cases, the total travel costs are minimized \eqref{E:objective}. Constraints \eqref{E:flow} are the flow conservation constraints. They ensure that each visited event is left again. Constraints \eqref{E:pick} take care that all pick-up locations are visited by a vehicle. Constraints \eqref{E:flow}--\eqref{E:pick} correspond to constraints \eqref{eq:degIn}--\eqref{eq:degOut} in the LB formulation. Due to constraint \eqref{E:K} at most $K$ vehicles leave the depot (corresponds to constraint \eqref{eq:depot}). Thus, at most $K$ tours are constructed. Constraints \eqref{E:time} focus on the time consistency. If a vehicle drives from location $v_1$ to location $w_1$ using the events $v$ and $w$, then the service at $w_1$ cannot start before the service at $v_1$ is finished ($B_v + s_{v_1}$) and the vehicle has driven to location $w_1$ ($t_{(v,w)}$; corresponds to constraints \eqref{eq:time}). Time consistency is considered separately for the first customer in a vehicle route (constraints \eqref{E:time1}). Constraints \eqref{E:tw1}--\eqref{E:tw3} focus on the time windows. While constraint \eqref{E:tw1} is the special case for the depot, events are separated according to their current location into constraints \eqref{E:tw2} (pick-up locations) and \eqref{E:tw3} (delivery locations). This is due to the fact that each location is represented by different events in the EB model. However, only one of the events corresponding to a location is used in a solution. See \citet{GKS22} for a detailed derivation of the time window constraints. Maximum ride time constraints are implemented analogously to the LB formulation, cf. constraints~\eqref{eq:mrt} and \eqref{E:mrt}. Note that $B_v \geq 0$ for all $v \in V$ because of constraints \eqref{E:tw1}--\eqref{E:tw3}, as $e_{v_1}\geq 0$ for all $v \in V$. 

Comparing the EB formulation with the LB formulation we find advantages for both of them (implicit formulation of precedence, pairing, and capacity constraints vs.\ fewer binary variables). We reformulate the EB model in the next sections to integrate the best of both models and show the advantages of the resulting formulations over the LB formulation theoretically afterwards.

\subsection{Location-augmented-event-based MILP Formulation} \label{IEBM}
In this section, we integrate the implementation of time consistency in the LB formulation into the EB formulation to obtain a more efficient formulation. 
\begin{subequations}
	\begin{align} 
		Z_{LAEB} = \min\;& \sum_{a \in A} c_a \cdot x_a \label{E2:objective}\\
		\text{s.t.}\;&\eqref{eq:B}, \eqref{eq:mrt}, \eqref{E:flow}-\eqref{E:K}, \text{ and } \eqref{E:binary}\quad&&\nonumber \\
		&\bar{B}_{j}\geq \bar{B}_{i}+ s_{i} + \bar{t}_{ij} - \bar{M}_{ij} \Biggl( 1 - \sum_{v:v_1=i,w:w_1=j} x_{(v,w)} \Biggr)\quad&&\forall i,j \in J\label{E2:time}
	\end{align}
\end{subequations}

Due to the construction of the event-based graph and
%Because of the graph formulation as well as
 constraints \eqref{E:pick} and \eqref{E:binary} there is exactly one node $w$ with $w_1 = i$ and $\sum_{a \in \delta^{in}(w)} x_a = 1$ while $\sum_{a \in \delta^{in}(w')} x_a = 0$ for all other nodes $w'$ with $w'_1 = i$. Therefore, only for the first one $M_{vw}(1 - x_{(v,w)})$ is 0 in constraints \eqref{E:time}. This constraint sets $B_w = \bar{B}_j$ while the others only ensure that feasible values are used for $B_{w'}$. This leads to the simplification of time consistency constraints \eqref{E2:time}. By this, time windows constraints \eqref{eq:B} and maximum ride time constraints \eqref{eq:mrt} can be ensured as in the LB formulation. As can be seen, the number of constraints in \eqref{eq:B}, \eqref{eq:mrt}, and \eqref{E2:time} is clearly smaller than in \eqref{E:time}--\eqref{E:mrt}. Moreover, the numerically unfavorable big-$M$-constraints are stronger in \eqref{E2:time} than in \eqref{E:time}, as
\[
\Biggl( 1 - \sum_{v:v_1=v'_1,w:w_1=w'_1} x_{(v,w)} \Biggr) \leq (1-x_{(v',w')}) \hspace{3cm} \forall (v',w') \in A
\]
while $\bar{M}_{ij} = M_{vw}$ with $v_1=i$ and $w_1=j$ with the classical choice
\[
\bar{M}_{ij} = \ell_{i} + s_{i} + \bar{t}_{ij} - e_{j}\hspace{2cm} \forall i,j \in \bar{J}
\]
and
\[
M_{vw} = \ell_{v_1} + s_{v_1} + t_{(v,w)} - e_{w_1}\hspace{2cm} \forall (v,w) \in A.
\]
Thus, the model formulation \eqref{eq:B}, \eqref{eq:mrt}, \eqref{E:flow}--\eqref{E:K}, \eqref{E:binary}, and \eqref{E2:objective}--\eqref{E2:time} is an improved version of \eqref{E:objective}--\eqref{E:binary}. In the next section, we present another location-augmented-event-based formulation with fewer binary variables.

\subsection{Aggregated location-augmented-event-based MILP Formulation} \label{SELF}
The LB formulation has the advantage that only $\mathcal{O}(n^2)$ binary variables are required while the EB formulation has $\mathcal{O}(n^{Q+1})$ binary variables (for $n \geq Q+1$). We use the equality
\begin{equation} \label{xbarx}
	\bar{x}_{ij} = \sum_{(v,w) \in A: v_1 = i \wedge w_1 = j} x_{(v,w)},
\end{equation}
which is already used in \eqref{E2:time}. Thereby, we get the following ALAEB formulation:
\begin{subequations} 
	\begin{align} 
		Z_{ALAEB} = \min \;& \sum_{i,j \in J} \bar{c}_{ij} \cdot \bar{x}_{ij}\label{E3:objective}\\
		\text{s.t.}\;&\eqref{eq:B}, \eqref{eq:mrt}, \eqref{eq:xbin}, \eqref{E:flow}-\eqref{E:K}, \text{ and } \eqref{xbarx}\quad&&\nonumber\\
		&\bar{B}_{j} \geq \bar{B}_{i} + s_{i} + \bar{t}_{ij} - \bar{M}_{ij} \left( 1 - \bar{x}_{ij} \right)\quad&&\forall i,j \in P \cup D\label{EL2:time} \\
		&0 \leq x_a \leq 1 \quad&&\forall a\in A\label{barxrel}
	\end{align}
\end{subequations}

While \eqref{E3:objective} corresponds to \eqref{eq:objFunc}, constraints \eqref{EL2:time} are equivalent to \eqref{E2:time} using \eqref{xbarx}. Thus, the difference between both location-augmented-event-based models is the fact that $\bar{x}$ variables are added as binary variables while $x$ variables are relaxed in constraints \eqref{barxrel}. For this reason, the model contains fewer binary variables than the LAEB formulation. Note, that given an integer solution $\bar{x}$ the tours are completely described. Thus, there are at most $K$ flows of size one from depot to depot through the event-based graph, i.e., $x_{a}$ variables are also integer. In total, the ALAEB formulation replaces an exponential number of constraints \eqref{eq:SEC} by an additional number of $\mathcal{O}(n^{Q+1})$ continuous variables.

\section{Theoretical Analysis} \label{sec:TA}
In this section, we investigate the properties of the LP relaxations of the LAEB and the ALAEB formulation in comparison to the LP relaxation of the LB formulation. We introduced the ALAEB formulation mainly to branch on a different set of variables.
%However, there is  no branching in the LP relaxation. Correspondingly, the first theorem shows that the LP relaxations of the LAEB and the ALAEB formulation are equivalent.
The following theorem shows, however, that the LP relaxations of the LAEB and the ALAEB formulation are equivalent.

\begin{Theorem} \label{equi}
The LP relaxations of the LAEB formulation \eqref{eq:B}, \eqref{eq:mrt}, \eqref{E:flow}--\eqref{E:K}, \eqref{E:binary}, and \eqref{E2:objective}--\eqref{E2:time} and of the ALAEB formulation \eqref{eq:B}, \eqref{eq:mrt}, \eqref{eq:xbin}, \eqref{E:flow}--\eqref{E:K}, \eqref{xbarx}, and \eqref{E3:objective}--\eqref{barxrel} are equivalent.
\end{Theorem}

The proof of Theorem \ref{equi} can be found in electronic companion \ref{PTequi}. We formulate and prove the following theorems only for the LAEB formulation, but due to Theorem~\ref{equi} the results hold also for the ALAEB formulation. The next theorem considers the special case that the time windows induce a unique order in time for every pair of locations.%all time windows reduce to a single point in time.

\begin{Theorem} \label{Net} 
	If the time windows $[e_i,l_i]$, $i \in \bar{J}$, fulfill the following conditions
\begin{enumerate}
	\item $l_{i^-} - e_{i^+} - s_i \leq L_i$ for all $i \in P$ and \label{C1}
	\item either $l_i + s_i + \bar{t}_{ij} \leq e_j$ or $e_i + s_i + \bar{t}_{ij} > l_j$ holds for all $i,j \in \bar{J}$, \label{C2}
\end{enumerate}
the LAEB formulation \eqref{eq:B}, \eqref{eq:mrt}, \eqref{E:flow}--\eqref{E:K}, \eqref{E:binary}, and \eqref{E2:objective}--\eqref{E2:time} is integral if all arcs $(v,w)$ with $v_1 = i$ and $w_1 = j$ are deleted if $l_j < e_i + s_i + \bar{t}_{ij}$.
%
%	If each customer location has a fixed time, that is, the time windows are limited to one unique point in time, i.e., $e_{j}= \bar{B}_{j} = \ell_{j}$ for all $j\in P \cup D$, the LAEB formulation \eqref{eq:B}, \eqref{eq:mrt}, \eqref{E:flow}--\eqref{E:K}, and \eqref{E:binary}, \eqref{E2:objective}--\eqref{E2:time} is integral if maximum ride time constraints are fulfilled and all arcs $(v,w)$ with $v_1 = i$ and $w_1 = j$ are deleted if $\bar{B}_{j} < \bar{B}_{i} + s_{i} + \bar{t}_{ij}$.
\end{Theorem}

The proof of Theorem \ref{Net} can be found in the electronic companion \ref{ProofNet}. Note that condition~\ref{C2}.\ holds in particular if each pick-up and delivery has a fixed time, that is, the time windows are limited to one unique point in time, i.e., $e_{j}= \bar{B}_{j} = \ell_{j}$ for all $j\in P \cup D$ (and time windows at depots are set appropriately). The theorem allows for several conclusions:
\begin{itemize}
	\item The sequencing of requests makes %Time windows make 
	the problem challenging. In fact, it is well-known that exact solution approaches for the DARP typically perform better for instances with tighter time windows which have a lower number of feasible sequences of requests.
	\item Condition \ref{C2}.\ ensures that $\bar{B}_{i} + s_{i} + \bar{t}_{ij} \leq \bar{B}_{j}$ holds for all arcs $(v,w)$ with $v_1 = i$ and $w_1 = j$ in the reduced graph. Besides the fact that constraints \eqref{E2:time} become unnecessary, this leads to a cycle-free graph structure. That is, there are no events any more which can be predecessors as well as successors of each other in different solutions. %We will work on the avoidance of subtours in Section~\ref{VS} as well.
	\item The example in the electronic companion \ref{ExlocMod} shows that Theorem \ref{Net} cannot be transferred to the LB formulation \eqref{eq:objFunc}--\eqref{eq:xbin}. This fact suggests that the LAEB and the ALAEB formulation are tighter than the LB formulation. 
\end{itemize}
The following theorem gives formal evidence for the last point.

\begin{Theorem} \label{T1} 
	Let $Z_{LB}^{rel}$ be the objective value of an optimal solution of the LP relaxation of the LB model \eqref{eq:objFunc}--\eqref{eq:xbin}, i.e., with $0 \leq \bar{x}_{ij} \leq 1$ for all $i,j \in \bar{J}$ instead of \eqref{eq:xbin}. Let further $Z_{LAEB}^{rel}$ be the objective value of an optimal solution of the LP relaxation of the LAEB model \eqref{eq:B}, \eqref{eq:mrt}, \eqref{E:flow}--\eqref{E:K}, \eqref{E:binary}, and \eqref{E2:objective}--\eqref{E2:time}, i.e., with $0 \leq x_a \leq 1$ for all $a \in A$ instead of \eqref{E:binary}. Then, $Z_{LAEB}^{rel} \geq Z_{LB}^{rel}$ holds. Moreover, there are instances in which $Z_{LAEB}^{rel} > Z_{LB}^{rel}$ holds.
\end{Theorem}

The proof of Theorem \ref{T1} can be found in the electronic companion \ref{PT1}. The proof highlights again that the LAEB and the ALAEB formulation implicitly ensure that pairing, precedence, and capacity constraints are fulfilled. The implicit implementation of capacity constraints leads even to a tighter polyhedron of the LP relaxation. Due to their exponential size, pairing and precedence constraints \eqref{eq:SEC} are typically not added upfront but only if they are violated. Therefore, the LAEB and the ALAEB LP relaxation are also tighter in this aspect. However, we still face a significant number of nodes and arcs in the event-based graph, which we compactify in Section~\ref{sec:BC}.

\section{Preprocessing and Branch-and-Cut Methods} \label{sec:BC}
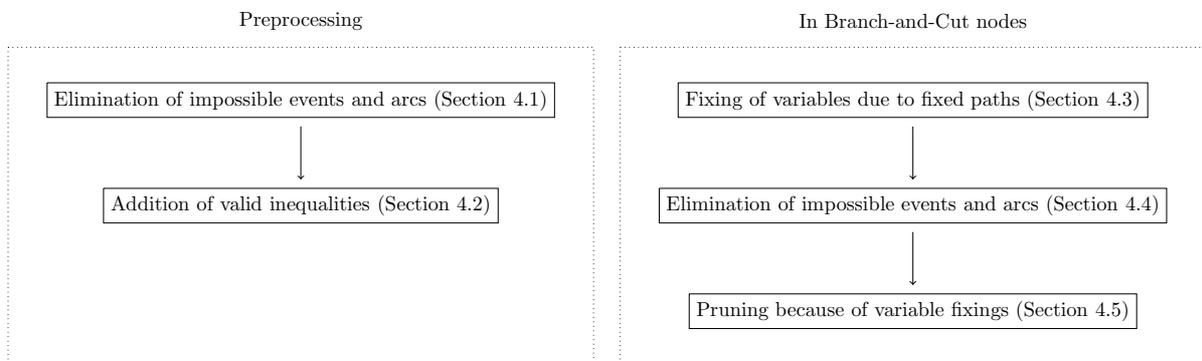
\begin{figure}[htbp]
	\centering
	\scalebox{0.7}{\begin{tikzpicture}
\node (A) at (0,6) {Elimination of impossible events and arcs (Section \ref{sec:Pre})};
\draw (A.north west) -- (A.north east) -- (A.south east) -- (A.south west) -- (A.north west);
\node at (0,4) (B) {Addition of valid inequalities (Section \ref{sec:VI})};
\draw (B.north west) -- (B.north east) -- (B.south east) -- (B.south west) -- (B.north west);

\node at (11.5,6) (E) {Fixing of variables due to fixed paths (Section \ref{sec:fixedpaths})};
\draw (E.north west) -- (E.north east) -- (E.south east) -- (E.south west) -- (E.north west);
\node at (11.5,4) (C) {Elimination of impossible events and arcs (Section \ref{sec:IBBC})};
\draw (C.north west) -- (C.north east) -- (C.south east) -- (C.south west) -- (C.north west);
%\node at (13,4) (D) {Improve local bounds on starting times in Branch-and-Cut nodes (Section 4.6)};
%\draw (D.north west) -- (D.north east) -- (D.south east) -- (D.south west) -- (D.north west);
\node at (11.5,2) (F) {Pruning because of variable fixings (Section \ref{sec:VFBC})};
\draw (F.north west) -- (F.north east) -- (F.south east) -- (F.south west) -- (F.north west);

\draw[->] (0,5.5) -- (0,4.5);
\draw[->] (11.5,5.5) -- (11.5,4.5);
\draw[->] (11.5,3.5) -- (11.5,2.5);

\draw[dotted] (-5.5,7) -- (5.5,7) -- (5.5,1) -- (-5.5,1) -- (-5.5,7);
\draw[dotted] (6,7) -- (17,7) -- (17,1) -- (6,1) -- (6,7);
\node at (0,7.5) {Preprocessing};
\node at (11.5,7.5) {In Branch-and-Cut nodes};
\end{tikzpicture}}
	\caption{Overview of developed methods}
	\label{fig:Flow_Chart}
\end{figure}

We present methods to reduce the size of the event-based graph in preprocessing and in branch-and-cut nodes as well as new types of additional valid inequalities in this section. Figure \ref{fig:Flow_Chart} gives an overview of the developed methods. 
%They can be divided into methods applied in preprocessing and in branch-and-cut nodes. 
First, we introduce preprocessing methods to eliminate impossible events and arcs in the event-based graph based on bounds for $B_v$ variables (Section~\ref{sec:Pre}). Moreover, we propose new valid inequalities to avoid subtours, incompatible events, and infeasible paths (Section~\ref{sec:VI}).
Second, methods are presented to improve the search in branch-and-cut nodes based on previous branching decisions, see Sections~\ref{sec:fixedpaths}--\ref{sec:VFBC}.

%Second, methods are introduced to improve the search in branch-and-cut nodes based on up to now branching decisions, \dg{e.g.\ fixing further arc variables to zero based on already fixed sequences of events (Section~\ref{sec:fixedpaths}), or eliminating nodes and arcs as a result of improving the bounds for $B_v$ variables derived in Section~\ref{sec:Pre} (Section~\ref{sec:IBBC}). One of the eliminated events might already been merged with other events. As this leads to an infeasible solution, the branch-and-cut node can be pruned (Section~\ref{sec:VFBC}.}

%The branching decisions lead to already fixed sequences of events (connected by arcs whose variables are fixed to 1). We use them to set further arc variables to zero if the corresponding arc cannot be part of the solution any more (Section~\ref{sec:fixedpaths}). Moreover, we can improve the bounds for $B_v$ variables derived in Section~\ref{sec:Pre} by including the branching decisions. Afterwards, these bounds are again used to eliminate impossible events and arcs (Section~\ref{sec:IBBC}). Finally, one of the eliminated events might already been merged with other events. As this leads to an infeasible solution, the branch-and-cut node can be pruned (Section~\ref{sec:VFBC}).

\subsection{Elimination of Impossible Event Nodes and Arcs in Preprocessing} \label{sec:Pre}
We already explained in Section~\ref{sec:ebMILP} that the number of events can be reduced by considering pairwise incompatibilities with respect to time windows and ride time constraints. We present a method to systematically compute earliest and latest starting times of service at all events. The earliest starting time $B_v^{LB}$ at node $v$ can be interpreted as the shortest path between depot event $(0, \ldots,0)$ and the considered event $(v_1, \ldots,v_Q)$, i.e., the fastest way to pick up first customers $v_2, \ldots, v_Q$ (for all $v_j \neq 0$) and $v_1$ afterwards if $v_1 \in P$ or deliver $v_1$ afterwards if $v_1 \in D$ while respecting the corresponding time windows. Analogously, the latest starting time of service $B_v^{UB}$ in $(v_1, \ldots,v_Q)$ is the latest time for start of service such that customers $v_1, \ldots, v_Q$, $v_j \neq 0$, can all be delivered within the time windows. Naturally, an event cannot be feasible if $B_v^{UB} < B_v^{LB}$. An event $v$ is also infeasible if there is a customer $i$ in the vehicle who cannot be delivered within the maximum ride time, i.e., if $\bar{B}_i^{UB} + s_i + L_i - s_{v_1} - t_{v_1,i^-} < B_v^{LB}$, whereat $\bar{B}_i^{UB} = \max_{v \in V_{i^+}} \{ B_v^{UB} \}$. Moreover, an arc $(v,w)$ is infeasible if $B_v^{LB} + s_{v_1} + t_{(v,w)} > B_w^{UB}$. Then, we can eliminate $(v,w)$ from the graph. 

Due to construction of the event-based graph,
\begin{equation} \label{BLB}
	B_v^{LB} = \max \left\{ e_{v_1}, \min_{w:(w,v) \in \delta^{in}(v) \wedge w \in \bigcup_{i=1}^n V_{i^+}} \left\{ B_w^{LB} + s_{w_1} + t_{(w,v)} \right\} \right\}
\end{equation}
and
\begin{align} \nonumber
	B_v^{UB} = &\min \left\{ \ell_{v_1} - \sum_{i \in R} \id_{\{ v_1 = i^-\}} \cdot \max\{0, \ell_{v_1} - (\bar{B}_i^{UB} + s_i + L_i)\} , \right.\\ \nonumber
	&\hspace{2cm}\max_{w:(v,w) \in \delta^{out}(v) \wedge w \in \bigcup_{i=1}^n V_{i^-}} \left\{ B_w^{UB} - t_{(v,w)} - s_{v_1} \right\}, \\ \label{BUB}
	&\hspace{4cm}\left. \min_{i \in R|\exists l \geq 2: v_l = i} \{ \bar{B}_i^{UB} + s_i+ L_i - s_{v_1} - t_{v_1,i^-} \} \right\}.
\end{align}
We can compute $B_v^{LB}$ and $B_v^{UB}$ in $\mathcal{O}(|V|)$ by evaluating the events $v$ systematically in the correct sequence to ensure that all predecessors (with $(w,v) \in A$) and successors (with $(v,w) \in A$), respectively, are evaluated before. $B_v^{LB}$ is computed in the following sequence:
\begin{enumerate}
	\item $v_1$ is a pick-up location and \ldots
	\begin{enumerate}
		\item $v_2 = \ldots = v_Q = 0$, i.e., all events where a customer is entering an empty vehicle.
		\item $v_2 \neq 0$, and $v_3 = \ldots = v_Q = 0$.
		\item $v_2, v_3 \neq 0$, and $v_4 = \ldots = v_Q = 0$.
		\item \ldots
		\item $v_2, \ldots, v_Q \neq 0$.
	\end{enumerate}
	\item $v_1$ is a delivery location.
\end{enumerate}
Note that the vehicle could also drive from a delivery location to a pick-up location. However, because of the triangle inequality this path cannot be shorter than not visiting neither the pick-up nor the delivery location of the corresponding customer. Therefore, it is reasonable to require $w \in \bigcup_{i=1}^n V_{i^+}$ in \eqref{BLB}. For the same reason we only need $B_w^{LB}$ of events where $w_1$ is a pick-up location to compute $B_v^{LB}$ for events where $v_1$ is a delivery location. 

Because of the triangle inequality the shortest path between an event $v$ and the end depot $(0, \ldots,0)$ is to deliver customers $v_1, \ldots, v_Q$ without picking up a new customer. When we first iterate over $v\in V$ to compute $B_v^{UB}$, maximum ride times cannot be considered yet, because we need the upper bounds $B_v^{UB}$ of all pick-up nodes to compute $\bar{B}_i^{UB}$. However, for the computation of $B_v^{UB}$, where $v$ is a pick-up node, we need the upper bounds of all delivery nodes. Hence, in the first loop, we omit maximum ride times from the computation of $B_v^{UB}$. The following sequence is most efficient to determine $B_v^{UB}$:
\begin{enumerate}
	\item $v_1$ is a delivery location and \ldots
	\begin{enumerate}
		\item $v_2 = \ldots = v_Q = 0$, i.e., all events where $v_1$ is the last customer leaving the vehicle.
		\item $v_2 \neq 0$, and $v_3 = \ldots = v_Q = 0$.
		\item $v_2, v_3 \neq 0$, and $v_4 = \ldots = v_Q = 0$.
		\item \ldots
		\item $v_2, \ldots, v_Q \neq 0$.
	\end{enumerate}
	\item $v_1$ is a pick-up location.
\end{enumerate}
In total, all events in $V$ have to be considered once to determine $B_v^{LB}$ and once to determine $B_v^{UB}$ without consideration of maximum ride times. After all upper bounds $B_v^{UB}$ have been computed, we are able to determine $\bar{B}_i^{UB}$, $i\in R$. We repeat the procedure to compute $B_v^{UB}$, this time including ride times. So, one iteration of the procedure requires an effort of $\mathcal{O}(\vert V\vert)$.

After computing $B_v^{UB}$ for an event $v$, we can directly check whether $B_v^{UB} < B_v^{LB}$ and delete the event if so. Due to our sequence to compute $B_v^{UB}$ values, the event is not considered in \eqref{BUB} any more to compute upper bounds for its predecessors. 

However, if an event $v$ with $v_1 \in P$ is deleted, lower bounds $B_w^{LB}$ of its successor events with $(v,w) \in A$ might also change. Therefore, we store all of these successors in a list and update their lower bounds afterwards. When going through the list, deleting an event $w$ can also lead to updated upper bounds for predecessors if $w_1 \in D$. Moreover, $\bar{B}^{UB}_i$ can change if $B_v^{UB}$ changes for a $v \in V_{i^+}$, $i\in R$. Thus, we need three lists, one for predecessors, one for successors, and one to store all $i$ for which $\bar{B}^{UB}_i$ has to be updated. Then, we run alternately through the lists in the sequence described above until the lists are empty. In the worst case, only one bound changes in each iteration.  

Whenever updating a lower or upper bound of event $v$, we also check feasibility of all in- and outgoing arcs ($(w,v) \in \delta^{in}(v)$ and $(v,w) \in \delta^{out}(v)$, respectively). In our tests, the computation of lower and upper bounds for all nodes and the detection of infeasible arcs was $<0.1s$ for instances with $n=100$ and $Q=6$.

We can use the bounds to add the following valid inequalities to the EB model formulation
\begin{align}
	&B_v^{LB} + (\ell_{i^-} - L_i - s_i - B_v^{LB}) \biggl(1 - \sum_{a\in\delta^{\text{in}}(v)} x_a \biggr) \leq B_v \leq B_v^{UB} \quad&&\forall i \in R, v \in V_{i^+}\label{bounds_eb_pickup} \\
	& B_v^{LB} \leq B_v \leq  e_{i^+}  + s_i +L_i +  (B_v^{UB} - (e_{i^+}  + s_i +L_i )) \sum_{a \in \delta^{in}(v)} x_a\quad&&\forall i \in R, v \in V_{i^-}\label{bounds_eb_delivery}
\end{align}
and the following valid inequalities to the location-augmented-event-based formulations:
\begin{equation} \label{bounds}
	\sum_{v:v_1=j} B_v^{LB} \cdot \sum_{a \in \delta^{in}(v)} x_a \leq \bar{B}_j \leq \sum_{v:v_1=j} B_v^{UB} \cdot \sum_{a \in \delta^{in}(v)} x_a \hspace{3cm} \forall j \in P \cup D
\end{equation}

\subsection{Addition of Valid Inequalities in Preprocessing} \label{sec:VI}
We can add several further valid inequalities to the model formulation identifying events which cannot occur simultaneously.
\subsubsection{Infeasible Paths}
In the literature, several authors introduced infeasible paths constraints \citep[e.g.][]{AFG00,Cor06,RCL07} based on sequences of request locations which cannot occur due to time windows. Certainly, if a path of locations $i \rightarrow j \rightarrow k$ would lead to a violated time window, this is also true for all event paths $v \rightarrow w \rightarrow u$ with $v_1 = i$, $w_1 = j$, and $u_1 = k$. However, it might be that the event path $v \rightarrow w \rightarrow u$ is infeasible although the location path $i \rightarrow j \rightarrow k$ is feasible (e.g. because in event $v$ another customer is in the vehicle who leads to a later departure time at $i$). Thus, it is possible to add a larger amount of infeasible path constraints based on events. 

Let $v$ be an event. Compute for all events $w$ with $(w,v) \in A$ the earliest starting time $B_{wv}^{LB} = B_w^{LB} + s_{w_1} + t_{(w,v)}$ at $v$ coming from $w$ and for all events $u$ with $(v,u) \in A$ the latest possible starting time of service $B_{vu}^{UB} = B_u^{UB} - s_{v_1} - t_{(v,u)}$ in $v$ to reach $u$ on time. If $B_{vu}^{UB} < B_{wv}^{LB}$, $x_{(w,v)} + x_{(v,u)} \leq 1$ is a valid inequality which can be lifted to
\begin{equation} \label{IP1}
	\sum_{w' \in V: (w',v) \in A \wedge B_{w'v}^{LB} \geq B_{wv}^{LB}} x_{(w',v)} + \sum_{u' \in V: (v,u') \in A \wedge B_{vu'}^{UB} \leq B_{vu}^{UB}} x_{(v,u')} \leq 1 \quad \forall v\in V\setminus V_0
\end{equation}
by including other predecessor (successor) events leading to a later earliest (earlier latest) starting time of service in $v$. 

In general, let $\bar{S} = \{(u^1,u^2), \ldots, (u^{\bar{m}-1},u^{\bar{m}})\}$ be a path $u^1 \rightarrow u^2 \rightarrow \cdots \rightarrow u^{\bar{m}}$, which is infeasible if $B_{u^1}^{LB} + \sum_{m = 1 }^{\bar{m}-1} s_{u^m} + t_{(u^m,u^{m+1})} > B_{u^{\bar{m}}}^{UB}$, i.e., if location $u^{\bar{m}}$ is not reached before its latest departure time. Then, $\sum_{a \in \bar{S}} x_a \leq \bar{m} - 2$ is a valid inequality. 

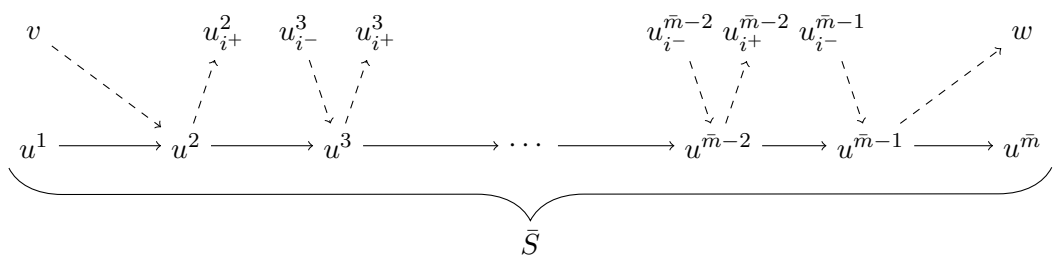
\begin{figure}[htbp]
	\centering
	\scalebox{1}{\usetikzlibrary{decorations.pathreplacing}
\begin{tikzpicture}
\node (A) at (0,0) {$u^1$};
\node (B) at (2,0) {$u^2$};
\node (B2) at (4,0) {$u^3$};
\node (C) at (6.5,0) {$\cdots$};
\node (D2) at (9,0) {$u^{\bar{m}-2}$};
\node (D) at (11,0) {$u^{\bar{m}-1}$};
\node (E) at (13,0) {$u^{\bar{m}}$};

\node (F) at (0,1.5) {$v$};
\node (G) at (2.5,1.5) {$u^2_{i^+}$};
\node (H) at (3.5,1.5) {$u^3_{i^-}$};
\node (G2) at (4.5,1.5) {$u^3_{i^+}$};

\node (G3) at (9.5,1.5) {$u^{\bar{m}-2}_{i^+}$};
\node (H3) at (8.5,1.5) {$u^{\bar{m}-2}_{i^-}$};
\node (H4) at (10.5,1.5) {$u^{\bar{m}-1}_{i^-}$};

\node (J) at (13,1.5) {$w$};

\draw[->] (A) to (B);
\draw[->] (B) to (B2);
\draw[->] (B2) to (C);
\draw[->] (C) to (D2);
\draw[->] (D2) to (D);
\draw[->] (D) to (E);
\draw[->,dashed] (F) to (B);
\draw[->,dashed] (B) to (G);
\draw[->,dashed] (H) to (B2);
\draw[->,dashed] (B2) to (G2);

\draw[->,dashed] (H3) to (D2);
\draw[->,dashed] (D2) to (G3);
\draw[->,dashed] (H4) to (D);

\draw[->,dashed] (D) to (J);

\draw[decorate,decoration={brace,raise=0mm,amplitude=20pt,mirror}] (A.south west) to node[below=7mm] {$\bar{S}$} (E.south east);

\end{tikzpicture}}
	\caption{Lifting for infeasible paths}
	\label{fig:InfeasiblePaths}
\end{figure}

This valid inequality can be lifted in two ways as Figure~\ref{fig:InfeasiblePaths} shows. In line with the argumentation above, $\sum_{a \in \bar{S}} x_a \leq \bar{m} - 2$ can be lifted by adding $x_{(v,u^2)}$ if $B_{vu^2}^{LB} \geq B_{u^1u^2}^{LB}$ and $x_{(u^{\bar{m}-1},w)}$ if $B_{u^{\bar{m}-1},w}^{UB} \leq B_{u^{\bar{m}-1}u^{\bar{m}}}^{UB}$ on the left side. Moreover, further paths between $u^2$ and $u^{\bar{m}-1}$ can be added if $(i^-, 0, \ldots, 0) \notin \{ u^2, \ldots, u^{\bar{m}-1} \}$ for all $i \in R$, i.e., the vehicle is not empty in between. In Figure~\ref{fig:InfeasiblePaths}, $u^m_{i^+}$ is the event where customer $i$ enters the vehicle and all customers who are in the vehicle after event $u^m$ are still there. Analogously, $u^m_{i^-}$ is the event where customer $i$ leaves the vehicle and all customers who are present in the vehicle directly before event $u^m$ are in the vehicle. After each of the path's locations the vehicle can leave it to visit the pick-up location of a further customer~$i$ and come back to the path when customer~$i$ is delivered. This might be directly, i.e., the path $u^m \rightarrow u^m_{i^+} \rightarrow u^{m+1}_{i^-} \rightarrow u^{m+1}$ replaces the path $u^m \rightarrow u^{m+1}$, or indirectly, i.e., the vehicle also visits further locations between pick-up and delivery of customer $i$ or drives from $i$'s delivery location to a later location of the original path. In the direct case, arcs $(u^m,u^m_{i^+})$ and $(u^{m+1}_{i^-},u^{m+1})$ replace arc $(u^m,u^{m+1})$ such that we have to weight both added arc variables with a factor $1/2$ to ensure feasibility of the lifted inequality. Note that we only ``replace'' the original arc conceptually in the solution, but we do not actually remove it from the lifted inequality. With respect to feasibility of the lifting, there are three cases:

\begin{itemize}
	\item In generalization of the already discussed situation, $u^m \rightarrow u^m_{i^+} \rightarrow \cdots \rightarrow u^{m+1}_{i^-} \rightarrow u^{m+1}$ replaces $u^m \rightarrow u^{m+1}$, which is a feasible lifting, as one arc with a factor of 1 is replaced by two arcs with a factor of $1/2$. Note that we do not need to consider the path between customer $i$'s pick-up and delivery. Due to the triangle inequality, the arrival time at $u^{m+1}$ cannot be earlier than in the original path independent of the path between customer $i$'s pick-up and delivery.
	\item $u^m \rightarrow u^m_{i^+} \rightarrow \cdots \rightarrow u^{m'}_{i^-} \rightarrow u^{m'}$ with $m' > m+1$ replaces the path $u^m \rightarrow \cdots \rightarrow u^{m'}$, which is a feasible lifting, as at least two arcs with a factor of 1 ($(u^m,u^{m+1})$ and $(u^{m'-1},u^{m'})$) are replaced by two arcs with a factor of $1/2$. Again, we do not need to consider the path between customer $i$'s pick-up and delivery here.
	\item If any path including arcs $(u^m_{i^-},u^m)$ and $(u^{m'},u^{m'}_{i^+})$ with $m' \geq m$ is part of the lifted path, the path is already infeasible due to precedence relations.
\end{itemize}
Note that the event $u^m_{i^+}$ ($u^m_{i^-}$) or arc $(u^m,u^m_{i^+})$ ($(u^m_{i^-},u^m)$) might not exist for some $i\in R$ and $m\in \{2, \ldots,\bar{m}-2\}$. Then, the corresponding arcs are simply not added to the valid inequality. As all other corresponding arcs are added, the lifting includes also paths where the delivery location visited right before the vehicle returns to the original path needs not to belong to the same customer as the pick-up location visited directly after leaving the original path. Together,
\begin{align} \nonumber
	&\sum_{(v,u^2) \in A, v\neq u_1: B_{vu^2}^{LB} \geq B_{u^1u^2}^{LB}} x_{(v,u^2)} + \sum_{a \in \bar{S}} x_a + \sum_{m=2}^{\bar{m}-2} \sum_{(u^m,w) \in A: w_1 \in P} \frac{1}{2} \cdot x_{(u^m,w)} \\ \label{IP2}
	&+ \sum_{m=3}^{\bar{m}-1} \sum_{(w,u^m) \in A: w_1 \in D} \frac{1}{2} \cdot x_{(w,u^m)} + \sum_{(u^{\bar{m}-1},w) \in A, w\neq u^{\bar{m}}: B_{u^{\bar{m}-1},w}^{UB} \leq B_{u^{\bar{m}-1}u^{\bar{m}}}^{UB}} x_{(u^{\bar{m}-1},w)} \leq \bar{m} - 2
\end{align}
is a valid inequality. We implemented inequalities \eqref{IP2} for the case of $\bar{m} = 4$. Details on the implementation of inequalities \eqref{IP1} and \eqref{IP2} can be found in electronic companion \ref{ImplIP}. 
In total, for each arc $(v,u)\in A$ we construct one inequality of type \eqref{IP1}, hence there are  $\mathcal{O}(|A|)$ valid inequalities. In the case of inequalities \eqref{IP2}, we construct each inequality from an arc $(u^2,u^3)\in A$ and $u^3$'s successors, hence there are $\mathcal{O}\bigl(|A| \cdot (|R| + |Q|)\bigr)$ valid inequalities.

\subsubsection{Vehicle Sharing} \label{VS}
Let $i$ and $j$ be two customers. If they use the same vehicle, only one of them can enter the vehicle while the other is already sitting in it. Both of these events are incompatible to an event where $i$ or $j$ leaves the vehicle before the other one entered it. Thus,
\begin{align} \nonumber
	&\sum_{(v,w) \in A: w_1 = j^+ \wedge \exists l: w_l = i} x_{(v,w)} + \sum_{(v,w) \in A: w_1 = i^+ \wedge \exists l: w_l = j} x_{(v,w)}  \\ \label{VS1}
	+& \sum_{(v,w) \in A: v_1 = i^- \wedge w_1 = j^+} x_{(v,w)} + \sum_{(v,w) \in A: v_1 = j^- \wedge w_1 = i^+} x_{(v,w)} 
	\leq 1 \;\;\; \forall i,j \in R: i<j
\end{align}
are valid inequalities. We have $\mathcal{O}(|R|^2)$ valid inequalities of type \eqref{VS1} and add all of them in the preprocessing. 

Analogously, only one of them can leave the vehicle while the other is still sitting in it or one of them can leave the vehicle before the other enters it. Hence,
\begin{align} \nonumber
	&\sum_{(v,w) \in A: w_1 = i^- \wedge \exists l: w_l = j} x_{(v,w)} +\sum_{(v,w) \in A: w_1 = j^- \wedge \exists l: w_l = i} x_{(v,w)} \\ \label{VS2}
	& +\sum_{(v,w) \in A: v_1 = i^- \wedge w_1 = j^+} x_{(v,w)} + \sum_{(v,w) \in A: v_1 = j^- \wedge w_1 = i^+} x_{(v,w)} \leq 1 \;\;\; \forall i,j \in R: i<j
\end{align}
are valid inequalities. We have $\mathcal{O}(|R|^2)$ valid inequalities of type \eqref{VS2} and add all of them in the preprocessing. 

Moreover, an event where $i$ or $j$ leaves the vehicle directly after the other one, is incompatible to an event where one of them left the vehicle directly before the other one entered it. Thus,
\begin{align} \nonumber
	&\sum_{(v,w) \in A: v_1 = i^- \wedge w_1 = j^-} x_{(v,w)} + \sum_{(v,w) \in A: v_1 = j^- \wedge w_1 = i^-} x_{(v,w)} \\ \label{VS3}
	+ &\sum_{(v,w) \in A: v_1 = i^- \wedge w_1 = j^+} x_{(v,w)} + \sum_{(v,w) \in A: v_1 = j^- \wedge w_1 = i^+} x_{(v,w)} 
	\leq 1 \;\;\; \forall i,j \in R: i<j
\end{align}
are valid inequalities. We have $\mathcal{O}(|R|^2)$ valid inequalities of type \eqref{VS3} and add all of them in the preprocessing. 

Furthermore, if we select an arc $(v,w)$ with $v_1 \in D$ and $w_1 \in P$, customers $v_1$ and $w_1$ cannot share a vehicle simultaneously. Thus,
\begin{align} \nonumber
	\sum_{(v,w) \in A: v_1 = j^- \wedge w_1 = i^+} x_{(v,w)} + &\sum_{(v,w) \in A: v_1 = i^- \wedge w_1 = j^+} x_{(v,w)} \\ \label{VS4}
	+&\sum_{(u',u) \in A: u_1 = k \wedge \exists l_1: u_{l_1} = i \wedge \exists l_2: u_{l_2} = j} x_{(u',u)} \leq 1
\end{align}
is a valid inequality for all pairs of customers $i,j \in R$, $i<j$, being pairwise compatible and a further location $k \in J \backslash \{ 0 \}$. We have $\mathcal{O}(|J|\cdot|R|^2)$ valid inequalities of type \eqref{VS4} and add all of them in the preprocessing.

\subsubsection{Customer Incompatibility} \label{VFITW}
With the pairwise incompatibilities presented in Section~\ref{sec:ebMILP} we can identify in advance customer pairs $i$ and $j$ which cannot be served by the same vehicle due to time windows. Our procedure in Section~\ref{sec:Pre} eliminates all events where $i$ and $j$ are simultaneously in the vehicle. If the paths $i^+ \rightarrow i^- \rightarrow j^+ \rightarrow j^-$ and $j^+\rightarrow j^- \rightarrow i^+ \rightarrow i^-$ are infeasible due to time windows, we can conclude that at most one of them can share the vehicle with another customer~$k$.% compatible with $i$ and $j$. 
Otherwise customers $i$ and $j$ would served by the same vehicle (not necessarily at the same time) due to transitivity. Thus,
\begin{align} \nonumber
	&\sum_{(u,v) \in A: v_1 = k^+ \wedge \exists l: v_l = i} x_{(u,v)} + \sum_{(u,v) \in A: v_1 = i^+ \wedge \exists l: v_l = k} x_{(u,v)}  \\ \label{CI1}
	&+ \sum_{(u,v) \in A: v_1 = k^+ \wedge \exists l: v_l = j} x_{(u,v)} + \sum_{(u,v) \in A: v_1 = j^+ \wedge \exists l: v_l = k} x_{(u,v)} \leq 1
\end{align}
is a valid inequality. We have $\mathcal{O}(|R|^3)$ valid inequalities of type \eqref{CI1}. However, the valid inequality is only relevant if $i$ and $j$ are incompatible, but both of them are compatible with $k$. We add all of them in the preprocessing.

\subsection{Fixing of Variables due to Fixed Paths}  \label{sec:fixedpaths}
Knowing which customers are together in a fixed path, i.e., a sequence of events connected by arcs whose variables are fixed to 1, also leads to further incompatible events if the depot is part of the path \citep{SP22}. We update the set of fixed paths when an upward branch, i.e., a branch where an $x$ or $\bar{x}$ variable is fixed to 1, is created. Let two fixed paths be given which start in the depot. As both have to be served by different vehicles, customers in one of the paths cannot share a vehicle with customers of the other path. Let customer $i$ be in the first and customer $j$ be in the second path. Then, we can eliminate all events $v$
% and adjacent arcs 
with $v_{l_1} \in \{i, i^+, i^-\}$ and $v_{l_2} \in \{j, j^+, j^-\}$ with $l_1,l_2 \in \{ 1, \ldots, Q \}$ and their incident arcs. The same is true if we consider two fixed paths ending in the depot. Thus, we check such incompatibilities whenever a fixed path is merged with another fixed path containing the depot.

\subsection{Elimination of Impossible Events and Arcs in Branch-and-Cut Nodes} \label{sec:IBBC}
Variable fixings in branch-and-cut nodes influence the lower and upper bounds computed in the preprocessing, see Section~\ref{sec:Pre}. While for the EB formulation and the LAEB formulation of the DARP, we branch on variables $x_a$, $a\in A$, using the ALAEB formulation branching takes place on variables $\bar{x}_{ij}$, $i,j\in \bar{J}$. If a variable $x_{(v,w)} $ is fixed to 0 or to 1, either directly because a new branch for the EB or the LAEB model is created, or as an indirect result of another fixed variable $x_{(v',w')}$ or $\bar{x}_{ij}$, the effects on upper and lower bounds are as follows: If $x_{(v,w)}$ is fixed to 0, $w$ cannot be the successor of $v$ any more. Thus, $w$ can be excluded in the maximum in \eqref{BUB}. Analogously, $w$ can be excluded in the minimum in \eqref{BLB} if $x_{(w,v)}$ is fixed to 0. If we fix $x_{(v,w)}$ to 1, $w$ is determined as the successor of $v$ such that we can replace \eqref{BLB} by
\begin{align} 
    B_w^{LB} = &\max \bigl\{ e_{w_1}, \; B_v^{LB} + s_{v_1} + t_{(v,w)} \bigr\} \tag{$6'$}\\
    \intertext{and \eqref{BUB} by}
	B_v^{UB} = &\min \left\{ \ell_{v_1} - \sum_{i \in R} \id_{\{ v_1 = i^-\}} \cdot 
      \max \bigl\{ 0,\; \ell_{v_1} - (\bar{B}_i^{UB} + s_i + L_i) \bigr\} ,\  B_w^{UB} - t_{(v,w)} - s_{v_1},\right. \nonumber\\ 
	&\hspace{5.5cm}\left. \min_{i \in R|\exists l \geq 2: v_l = i} \bigl\{ \bar{B}_i^{UB} + s_i+ L_i - s_{v_1} -t_{v_1,i^-} \bigr\} \right\} \tag{$7'$}.
\end{align}
In other words, the events $v$ and $w$ are merged to a single event if $x_{(v,w)} = 1$ \citep[cf.][]{SP22}.
Moreover, we eliminate events and arcs from the graph which have become infeasible due to variable fixings in a branch-and-cut node. For example, if arc $x_{(v,w)}$ has been fixed to 1 and $w$ is a pick-up node, we know that all other events $u\in V$, $u_1 = w_1$, $u\neq w$, are infeasible for the subtree rooted in the current branch-and-cut node. 
%\todo[inline]{Kommentar für Michael: was zu strong branching/ Zuordnungsproblem schreiben: wenn Kante $x_{(v,w)}$ auf 1 fixiert, dann können alle anderen Knoten in $V_{v_1}$ und $V_{w_1}$ entfernt werden (z.B. hier einbauen?)}
%\subsection{Improving Bounds on Starting Times in Branch-and-Cut Nodes} \label{sec:improvebounds}
%After we ran all procedures described in Sections~\ref{sec:IBBC}--\ref{sec:fixedpaths}
After we ran the procedures described in this and the previous section, we store all events $v$ for which either $B_v^{LB}$ or $B_v^{UB}$ changed and update lower and upper bounds for the beginning of service of all possible predecessor and successor events as described at the end of Section~\ref{sec:Pre}. Moreover, we check feasibility of events $v$ and arcs $(v,w)$, i.e., check if $B_v^{LB} \leq B_v^{UB}$ and $B_v^{LB} + s_{v_1} + t_{(v,w)} \leq B_w^{UB}$, respectively. 

If fixing of variables $x_{(v,w)}$ or the update of lower and upper bound results in an infeasible event, the event is removed from the graph together with its incident arcs.
% together with its incident arcs. If the update results in an infeasible arc, the arc is removed. 
Finally, we add new inequalities of types \eqref{bounds_eb_pickup}--\eqref{bounds} if bounds changed, and inequalities $x_a \leq 0$ if the arc $a$ has been declared infeasible.

\subsection{Pruning because of Variable Fixing in Branch-and-Cut Nodes} \label{sec:VFBC}
If one of the eliminated events $v$ is already merged with other events, i.e., is incident to an arc for which the corresponding arc variable is fixed to 1, the branch-and-cut node can be pruned.

\section{Computational Study} \label{sec:CS}
In this section, we evaluate the efficiency of the presented MILP models, lower and upper bounding strategies, and valid inequalities. 
Comparing the EB and the LAEB model, the same branching decisions can be applied, as in both models branching takes place on variables $x_a$, $a\in A$. The LAEB model has fewer variables, as we replace variables $B_v$, $v\in V$, by variables $\bar{B}_j$, $j \in \bar{J}$. Hence, we expect the LAEB model to be more efficient.
Comparing our two location-augmented-event-based formulations, the LAEB model has fewer variables but more binary variables. Moreover, variable fixing in Section~\ref{sec:VFBC} can more likely be applied, as we directly branch on the arcs of the event-based graph. Therefore, a branch where an $x_a$, $a \in A$, variable is fixed to 1 directly implies that both adjacent events occur. This branching step contains significantly more information than just fixing a single relation between two locations. On the other hand, a sequence of variables $\bar{x}_{ij}$ has to be fixed to 1 to ensure that an event occurs if we branch on $\bar{x}_{ij}$ variables, but we have significant fewer $\bar{x}_{ij}$ variables. Taking everything together, it is unclear which advantage predominates. We compare all three formulations in a computational study. 

Our computations were performed on an Intel Core i7-8700 CPU, 3.20\,GHz, 32\,GB memory and implemented in C++ using CPLEX 12.10. Note, that the execution of CPLEX was limited to one thread only and the search method was limited to traditional branch-and-cut.% due to the intervention in CPLEX's branch-and-cut algorithm. 
%The code was written.
The time limit for the solution in all tests was set to $7200$ seconds. Throughout this section, we use some abbreviations which are summarized in Table~\ref{tab:abbrv}. We first analyze the performance of the presented MILP formulations (Section~\ref{sec:com_milp}). Then, our preprocessing components are validated in Section~\ref{sec:test_preprocessing}. Afterwards, we test the components for the branch-and-cut algorithm (Section~\ref{sec:bcalg}). We present our final results on benchmark instances in Section~\ref{sec:resultsbench}.

\subsection{Comparing the MILP formulations}
\label{sec:com_milp}
The EB, LAEB, and ALAEB formulations are compared against each other using the a- and b-benchmark instances introduced in \cite{Cor06} and \cite{RCL07}. Furthermore, we use an extended version of these instances denoted by a-X and b-X \citep{GI15}, where the time window length is doubled by postponing the upper bound of the pick-up or delivery time window by 15~minutes, i.e., $\ell_i \mathrel{+}=15$. Due to the larger time windows these instances are harder to solve. 
%The extended instances a-X and b-X were introduced in \citep{GI15}. 
Note, that we do not increase the vehicle capacity as done in \cite{GI15}.  The instances are denoted as a$K$-$n$, b$K$-$n$ or a$K$-$n$-X and b$K$-$n$-X, respectively, where $K$ indicates the number of vehicles and $n$ denotes the number of requests. In the a- and a-X-instances, $Q=3$ and $L_i = 30$ and $q_i=1$ for all $i\in R$, whereas in the b- and b-X-instances $Q=6$ and $L_i = 45$ and $q_i\in\{1,\ldots,6\}$ for all $i\in R$.  

The results are presented in Tables \ref{tab:comparison_ab} and \ref{tab:comparison_aXbX}.  On the smaller as well as on the larger benchmark instances, the LAEB model outperforms the other two formulations, but its superiority becomes most evident when comparing its performance to the other formulations on the harder a-X and b-X instances. 
The average root node gap for all three formulations is very small and ranges from 1.1\% to 5.4\%, which demonstrates that, although some instances still take a long time to solve, the MILP formulations are already very tight. To sum it up, it seems that the smaller number of integer variables in the ALAEB formulation does not make up for the loss of information compared to the EB or LAEB formulation when branching on $\bar{x}_{ij}$ variables instead of $x_{(v,w)}$. Comparing the EB and the LAEB model, the replacement of variables $B_v$, $v\in V$, (of which there may be $\mathcal{O}(n^Q)$  \citep[see][]{GKS22}) by variables $B_j$, $j  \in J$, (of which there are $2\, n$ variables) as well as the tighter big-M formulation explain the presented speedup. 
Since the LAEB formulation improves on the EB formulation and turns out to be the best formulation, we restrict ourselves to this formulation in the remaining parts of the numerical tests. 

\subsection{Validating Preprocessing Components}\label{sec:test_preprocessing}
We proposed a preprocessing procedure to eliminate infeasible nodes and arcs (Section~\ref{sec:Pre}) and three types of new valid inequalities (Section~\ref{sec:VI}). 
%To analyze the individual and joint effect of the preprocessing components we conducted numerical tests with all possible subsets of the components. The test set consists of the following instances:
To evaluate the effect of the aforementioned methods, we use the following test set:
\begin{itemize}
	\item instances a6-72, a8-80, b6-72, and b8-96 introduced above
	\item instances Q3n80\_2, Q3n80\_3, Q3n80\_4, Q6n80\_5, Q6n100\_1, and  Q6n100\_5 from the artificial instance test set used in \cite{GKS22}
	\item five real data instances from HolMichApp containing one full day of completed trips each 
\end{itemize}
The characteristics of the instances are summarized in Table~\ref{tab:training_instances}. The instances have been selected to provide adequate diversity in terms of share of requests per vehicle, share of number of requests per service duration, and length of time windows. 

First, we test the efficiency of graph preprocessing (GP). On average, the number of nodes is reduced by 32\%, and the number of arcs by 12\%. As the number of nodes and arcs translates directly to the number of variables in the MILP, GP leads to a reduction of 22\% of the variables in total.  Table~\ref{tab:pre} shows the results for a combination of GP with different subsets of valid inequalities. The first row shows the total solution time for all instances when no preprocessing components are switched on. 
In the following rows, we test all components and all but one component to be switched on and compare computation times to the first row, where no preprocessing takes place. We achieve an average deviation in computation time ranging from -17 to -52\% over all test instances and an average deviation ranging from -64\% to -82\% over the test instances where computation times are improved. The significant difference between the average deviation over all instances and the average deviation only over the instances with an improvement underlines the heterogeneity of the test set. 
%In the next row, all components are switched on leading to an average reduction of computational time of 17\% over the test instances and an average reduction of 64\% over the 12 test instances, where computation times are compared to the first row without any preprocessing methods. The significant difference between the average deviation over all instances and the average deviation only over the instances with an improvement underlines the heterogeneity of the test set. 
% Since it is not possible to switch on infeasible path constraints \eqref{IP1} and \eqref{IP2} without switching on graph preprocessing (GP) due to the lower and upper bounds needed in these inequalities, the third row contains three zeros. 
%We were able to achieve an average improvement in computation time ranging from 21 to 52\% with these settings, only for the case where all valid inequalities but IP1 and IP2 were switched on and GP was switched off, there was an overall increase in computational time. Since this is the only row where GP was switched off, the last row in this part of the table shows results for the case that all valid inequalities but IP1 and IP2 were switched on and GP is switched on, leading to an average improvement of 10\% in computational time. 
%From the first part of the table, we conclude that GP has a positive influence on computational times. 
In the next part of the table, we compare the influence of only one of the valid inequalities switched on (and GP switched on, since otherwise the comparison would be unfair, as IP1 and IP2 cannot be switched on without GP). 
%We achieved average reductions in computational time ranging from 4\% to 58\%. Only the combination of GP and VS3 led to an average increase of 4\%, but on the other hand, some of the rows above including VS3 led to the strongest improvement in computational time. Thus, all of the introduced preprocessing steps have a benefit at least for a subset of instances. 
The improvements in the second part of the table are not as strong as the improvement in the first part, indicating that a combination of multiple valid inequalities and GP is most promising. Taking a look at the average root node gap, we are able to reduce it from 2.9\% to 1.9\% with multiple combinations of preprocessing steps. This again highlights the tightness of the LAEB formulation, but also underlines the efficiency of the preprocessing methods. It turns out that for this set of test instances, the best combination is the one given in the bottom row of the table, leading to an average reduction of 58\% in CPU. Although the test instances are very diverse, we cannot exclude that for other sets of test instances another combination prevails, so that we suggest to use a combination of all of the presented preprocessing methods. 

After verifying the new presented MILP formulation and preprocessing methods, we next examine the influence of the presented strategies in branch-and-cut nodes. 

\subsection{Testing Branch-and-Cut Algorithm} \label{sec:bcalg}
The lower part of Table~\ref{tab:abbrv} gives an overview of the methods introduced in Sections~\ref{sec:fixedpaths}--\ref{sec:VFBC} to reduce the size of the graph in branch-and-cut nodes. Following our strategy to evaluate individual components, we solved the LAEB model together with all of the presented preprocessing steps once for each subset of branch-and-cut components on the set of test instances introduced in the previous section. The results are listed in Table~\ref{tab:usercuts}.

The first row of Table~\ref{tab:usercuts} shows the total computational time for the LAEB model and preprocessing steps. In the next rows, all three, two or only one of the presented branching techniques are switched on. Note, that we did not include the case where neither LBsUBs or FA are switched on, since in this case no valid inequalities would be added. We were able to achieve further improvements of about 25\% compared to the total computational time after applying preprocessing techniques in about a third of the instances. However, taking the whole set of test instances into account, there is an overall increase ranging from 54\% to 66\% in CPU. This may be explained by the fact that the overhead caused by the implementation of our cuts in CPLEX is not balanced by the profit from using additional cuts. Also the average root node gap was very small already such that we might not get low enough in the branch-and-cut tree to fully use the benefit of the presented methods. Nevertheless, the results underline again the heterogeneity of the instances and show that the presented methods are very effective on certain instances. 

\subsection{Results on benchmark instances} \label{sec:resultsbench}
In this section, we show results on the benchmark instances using the LAEB formulation and the preprocessing methods. The results are shown in Table~\ref{tab:1thread_ab_aXbX}. On the a-instances, the total computational time was reduced by 94\% compared to using the plain LAEB or the plain EB formulation. On the b-instances, the computational time was reduced by 59\% and 79\%, respectively (compare Table~\ref{tab:comparison_ab}). On the a-X instances, the computational time was reduced by 27\% compared to using the plain LAEB formulation and by 41\% compared to the EB formulation. On the b-X instances, the computational time was reduced by 92\% and 96\%, respectively. These results prove the efficiency of the new LAEB formulation and the high impact of the proposed preprocessing techniques. Although computation times on different computers are only comparable to a limited extent, our results are highly competitive with the branch-and-cut algorithm by \cite{GI15}:  While some of the a-instances were solved faster by our approach and some of them were solved faster by their branch-and-cut algorithm, the solution times on the b-instances are smaller by about 50 orders of magnitude. However, our solution times cannot compete with the speed of the branch-and-cut algorithm developed by \cite{RF21} but could still prove useful in praxis, as the implementation of a MILP formulation and preprocessing methods is very fast and easy compared to the implementation of a branch-and-cut algorithm. Moreover, the results show that our methods are able to solve instances of substantial size within seconds. Thus, they can be used in a practical dynamic setting of a ridepooling provider within a rolling horizon approach like it is shown in \cite{GKS21}. This is especially true, as our results show that the model performs even better for instances where groups are transported (b and b-X instances) which is an important case for ridepooling providers. 
The results in this section are still for one thread only. Since we do not use user-defined cuts if we omit the methods presented in Sections \ref{sec:fixedpaths}--\ref{sec:VFBC}, we do not need to limit the number of threads used by CPLEX or to restrict the search method to traditional branch-and-cut. Using 12 threads the whole set of a- and b-instances was solved in 62 seconds. The corresponding results for the a-X and b-X instances show that medium-sized and some of the large instances of the harder benchmark set can be solved within a few seconds (compare Table~\ref{tab:12threads_aXbX} in the electronic companion).

\section{Conclusion} \label{sec:Con}
The new presented MILP formulations combine the existing state-of-the-art formulations, the location- and the event-based formulation of the DARP. As shown in Theorem~\ref{T1}, the LAEB formulation is tighter than the LB formulation. Computational tests on large benchmark instances show that, using the LAEB formulation, computational times can be reduced by 19\% (a-X instances) and 47\% (b-X instances), respectively, compared to the EB formulation. Additionally, graph preprocessing and the introduction of new valid inequalities which eliminate subtours, infeasible events, and infeasible paths strongly further improve computational times by 27\% (a-X) instances and 92\% (b-X), respectively. An average root node gap of 1.6\% proves that our formulation and preprocessing methods generate a very tight MILP model. Our methods are useful for application in praxis, as the implementation is fast and easy and even medium-sized instances of the harder benchmark set can be solved within a few seconds (compare Table~\ref{tab:12threads_aXbX} in the electronic companion).

While the literature mostly focuses on the static DARP, on-demand ridepooling services need to include customer requests on time when they arise. \citet{GKS21} showed that the dynamic DARP can be solved using a rolling-horizon algorithm in which the EB formulation is updated and solved whenever new requests arrive. Since the LAEB formulation improves on the EB formulation, this strategy could be adapted using the LAEB formulation as well.

As shown in Theorem~\ref{Net} the LAEB formulation is integral if time windows fulfill additional conditions (including the special case if they reduce to a single point in time). \citet{BJM19} used an analogue result for a taxi setting to solve the dynamic problem setting. A direction for future research is to use our result in Theorem~\ref{Net} and the ideas by \citet{BJM19} for the dynamic variant of the DARP. Moreover, instances with single point time windows lead to cycle-free event graphs. This special structure might also be an interesting direction for future research, as there are practical settings like the case where all customers have the same pick-up or delivery location as in the transition to public transport with fewer cycles in the event graph. Future research could evaluate how this influences the solution performance.

\paragraph{Acknowledgement}This work was partially supported by the state of North Rhine Westphalia (Germany) within the project "bergisch.smart\_mobility".

\bibliographystyle{apa}
\bibliography{Event_based}

%%%%%%%%%%

\newpage
\section*{Tables}

\setlength{\tabcolsep}{8pt}
\begin{table}[ht]
	\centering\footnotesize
	\begin{threeparttable}		
		\begin{tabular}{lp{12cm}}
			\toprule
			Abbreviation & Explanation/ Reference\\
			\midrule
			\multicolumn{2}{l}{\textbf{General abbreviations}} \\
			Inst. & Name of test instance \\
			Obj. & Objective Value \\
			EB & Event-based MILP \\
			LAEB & Location-augmented-event-based MILP \\
			ALAEB & Aggregated location-augmented-event-based MILP \\
			N/A & Not applicable (no integer solution found within the timelimit) \\
			\# Improvements & Number of instances with improved computational time \\
			CPU & CPU time in seconds \\
			Avg. Deviation & Average deviation of CPU time compared to the benchmark (first row of the corresponding table) \\
			Avg.* Deviation & Avg. Deviation only on instances with improved computational time \\
			rGap & root node gap \\
			Avg. rGap & Average root node gap \\[2ex]
			\multicolumn{2}{l}{\textbf{Preprocessing}} \\
			GP & Graph preprocessing, see Section~\ref{sec:Pre} \\
			VS1--VS4 & Vehicle Sharing, see \eqref{VS1}--\eqref{VS4} \\ 
			CI1 & Customer Incompatibility, see \eqref{CI1} \\
			IP1--IP2 & Infeasible Paths, see \eqref{IP1}--\eqref{IP2} \\[2ex]
			\multicolumn{2}{l}{\textbf{Branch-and-Cut Algorithm}} \\
			LBsUBs & add new inequalities of type  \eqref{bounds_eb_pickup}--\eqref{bounds} if bounds changed,  see Section~\ref{sec:IBBC} \\[2ex]
			FA & add inequalities $x_a\leq 0$ for infeasible arcs $a$, see Section~\ref{sec:IBBC} \\[2ex]
			FP & update the set of fixed paths and eliminate events (and incident arcs) which represent users of two different paths connected to the start or end depot  together in the vehicle, see Section~\ref{sec:fixedpaths} \\
			\bottomrule
		\end{tabular}
		\caption{Abbreviations used in Section~\ref{sec:CS}.}
		\label{tab:abbrv}
	\end{threeparttable}
\end{table} 

\setlength{\tabcolsep}{8pt}
\begin{table}[ht]
	\centering\footnotesize
	\begin{threeparttable}
		\begin{tabular}{lcccccc}
			\toprule
			Inst. & $n$ & $K$ & $T$ & $TW$ & $L_i$ & $Q$ \\						 
			\midrule
			a6-72 & 72 & 6 & 720 & 15 & 30 & 3\\
			a8-80 & 80 & 8 & 600 & 15 & 30 & 3\\
			b6-72 & 72 & 6 & 720 & 15 & 45 & 6\\
			b8-96 & 96 & 8 & 720 & 15 & 45 & 6 \\
			\midrule
			Q3n80\_2 & 80 & 11 & 240 & 15 & 1.5 $\bar{t}_i$ & 3 \\
			Q3n80\_3 & 80 & 11 & 240 & 15 & 1.5 $\bar{t}_i$ & 3 \\
			Q3n80\_4 & 80 & 10 & 240 & 15 & 1.5 $\bar{t}_i$ & 3 \\
			Q6n80\_5 & 80 & 12 & 240 & 15 & 1.5$\bar{t}_i$  & 6 \\
			Q6n100\_1 & 100 & 17 & 240 & 15 & 1.5$\bar{t}_i$ & 6 \\
			Q6n100\_5 & 100 & 14 & 240 & 15 & 1.5$\bar{t}_i$ & 6 \\
			\midrule
			2021-10-05 & 85 & 8 & 960 & 25 & $\bar{t}_i + \max(10, 0.75\bar{t}_i)$ & 6 \\
			2021-10-07 & 89 & 8 & 960 & 25 & $\bar{t}_i + \max(10, 0.75\bar{t}_i)$ & 6 \\
			2021-10-09 & 97 & 12 & 1080 & 25 & $\bar{t}_i + \max(10, 0.75\bar{t}_i)$ & 6 \\
			2021-10-16 & 120 & 12 & 1080 & 25 & $\bar{t}_i + \max(10, 0.75\bar{t}_i)$ & 6 \\
			2021-10-18 & 68 & 2 & 960 & 25 & $\bar{t}_i + \max(10, 0.75\bar{t}_i)$ & 6 \\
			\bottomrule		
		\end{tabular}
		\caption{Characteristics of the test instances.}
		\label{tab:training_instances}
	\end{threeparttable}
\end{table}

\setlength{\tabcolsep}{3pt}
\sisetup{table-align-text-post = false, table-number-alignment=center,table-format = 4.1}
\begin{table}[ht]
	\centering\footnotesize
	\begin{threeparttable}		
		\begin{tabular}{cSS[table-format = 1.1]SS[table-format = 1.1]S[table-format = 1.1]SS[table-format = 1.1]S[table-format = 1.1]SS[table-format = 1.1]}
			\toprule
			& & \multicolumn{3}{c}{EB} & \multicolumn{3}{c}{LAEB} & \multicolumn{3}{c}{ALAEB} \\
			\cmidrule(lr){3-5} \cmidrule(lr){6-8} \cmidrule(lr){9-11}
			$\text{Inst.}$ &  $\text{Obj.}$  & $\text{Gap}$  & $\text{CPU}$ & $\text{rGap}$ &   $\text{Gap}$  & $\text{CPU}$ & $\text{rGap}$ &   $\text{Gap}$  & $\text{CPU}$ & $\text{rGap}$\\
			\midrule
			a2-16 & 294,3 &&0,03&0,8\%  &&0,01&0,8 \% &&0,03&0,8 \% \\
			a2-20 & 344,9 &&0,03&0,0\%  &&0,01&0,0\%  &&0,02&0,0 \%  \\
			a2-24 & 431,1 &&0,09&1,9\%  &&0,03&2,1\%  &&0,06&2,1  \% \\
			a3-18 & 300,5 &&0,04&1,8\%  &&0,06&1,8\%  &&0,04&1,8  \% \\
			a3-24 & 344,8 &&0,14&0,5\%  &&0,08&2,8\%  &&0,3&3,0 \%  \\
			a3-30 & 494,8 &&0,05&0,0\%  &&0,02&0,0\%  &&0,08&0,6 \%  \\
			a3-36 & 583,2 &&0,21&2,5\%  &&0,13&2,5\%  &&0,16&2,5  \% \\
			a4-16 & 282,7 &&0,35&0,7\%  &&0,04&0,7\%  &&0,05&0,7 \%  \\
			a4-24 & 375,0 &&0,05&0,0\%  &&0,02&0,0\%  &&0,03&0,1 \%  \\
			a4-32 & 485,5 &&0,43&1,8\%  &&0,12&2,0\%  &&0,24&2,1  \% \\
			a4-40 & 557,7 &&0,83&1,7\%  &&1,11&1,8 \% &&2,04&2,8  \% \\
			a4-48 & 668,8 &&0,39&2,0\%  &&0,47&2,2\%  &&2,66&4,1 \%  \\
			a5-40 & 498,4 &&0,24&1,1\%  &&0,11&1,2\%  &&0,32&1,2 \%  \\
			a5-50 & 686,6 &&1,85&2,6\%  &&1,8&2,6\%  &&77,42&4,0 \%  \\
			a5-60 & 808,3 &&0,72&2,0\%  &&0,93&2,1\%  &&16,71&2,2 \%  \\
			a6-48 & 604,1 &&0,68&0,7\%  &&0,42&0,6 \% &&3,85&2,2 \%  \\
			a6-60 & 819,3 &&9,61&2,4\%  &&9,59&2,4\%  &&85,33&3,2 \%  \\
			a6-72 & 916,1 &&16,75&2,5\%  &&28,49&2,6 \% &&656,72&3,5 \% \\ 
			a7-56 & 724,0 &&4,85&1,3\%  &&1,76&1,6 \% &&54,59&3,1 \%  \\
			a7-70 & 875,7 &&14,99&1,2\%  &&4,47&1,3 \% &&199,76&2,2 \%  \\
			a7-84 & 1033,3 &&20,72&1,9\%  &&16,59&2,1 \% &&1535,79&4,6  \%\\ 
			a8-64 & 747,5 &&8,35&2,0\%  &&5,59&1,9\%  &&420,34&2,9 \%  \\
			a8-80 & 945,8 &&23,06&2,7\%  &&15,84&2,8 \% &1,0\% &2h&4,5 \% \\ 
			a8-96 & 1229,7 &1,0\% &2h&4,1\%  &1,1\% &2h&4,1\%  &2,9 \%&2h&5,2 \%  \\
			\cmidrule(lr){1-1}	\cmidrule(lr){4-4}	\cmidrule(lr){5-5}\cmidrule(lr){7-7}\cmidrule(lr){8-8}	\cmidrule(lr){10-10}\cmidrule(lr){11-11}
			Total/Avg. &             && \text{7304} &1,6\% & &\text{7287} & 1,8\%& &\text{17456} & 2,5\%\\
			\cdashlinelr{1-11}\addlinespace[2mm]
			b2-16&309,4&&0,03&0,6\%  &&0,07&0,6\%  &&0,04&0,6\%    \\
			b2-20&332,7&&0,01&0,0\%  &&0,04&0,0\%  &&0,01&0,0\%    \\
			b2-24&444,7&&0,05&1,7\%  &&0,05&1,8\%  &&0,05&1,8 \%   \\  
			b3-18&301,6&&0,04&2,0\%  &&0,05&2,2\%  &&0,04&2,2 \%   \\
			b3-24&394,5&&0,26&2,2\%  &&0,19&2,2\%  &&0,16&2,2\%    \\
			b3-30&531,4&&0,03&0,5\%  &&0,03&0,5\%  &&0,03&0,5 \%   \\
			b3-36&603,8&&0,05&0,3\%  &&0,04&0,3\%  &&0,04&0,3 \%   \\
			b4-16&296,9&&0,01&1,6\%  &&0,03&1,7\%  &&0,01&1,7 \%   \\
			b4-24&371,4&&0,04&0,7\%  &&0,04&0,8\%  &&0,05&0,8\%    \\
			b4-32&494,9&&0,03&0,0\%  &&0,03&0,0\%  &&0,03&0,0 \%   \\
			b4-40&656,6&&0,09&0,3\%  &&0,07&0,3\%  &&0,13&0,3\%    \\
			b4-48&673,8&&0,76&0,9\%  &&0,42&1,0\%  &&2,42&1,1\%    \\
			b5-40&613,7&&0,17&0,4\%  &&0,09&0,5\%  &&0,1&0,4 \%   \\
			b5-50&761,4&&0,42&1,4\%  &&0,2&1,5\%  &&0,55&1,6 \%   \\
			b5-60&902,0&&0,99&1,8\%  &&2,16&1,9\%  &&7,74&1,9\%    \\
			b6-48&714,8&&0,22&0,5\%  &&0,12&0,5\%  &&0,13&0,5\%    \\
			b6-60&860,0&&0,28&0,6\%  &&0,17&0,6\%  &&0,25&0,6 \%   \\
			b6-72&978,5&&8,85&1,0\%  &&5,64&1,0\%  &&26,95&1,0\%    \\
			b7-56&824,0&&6,79&0,9\%  &&4,01&0,9\%  &&22,72&0,9 \%   \\
			b7-70&912,6&&2,51&2,0\%  &&2,77&2,3\%  &&9,12&2,5 \%   \\
			b7-84&1203,3&&2,87&1,2\%  &&3,07&1,2\%  &&6,3&1,3 \%   \\
			b8-64&839,9&&2,27&2,3\%  &&1,1&2,3 \% &&5,56&2,3\%    \\
			b8-80&1036,4&&3,49&1,3\%  &&1,46&1,9\%  &&16,89&1,9 \%  \\ 
			b8-96&1185,6&&78,05&1,4\%  &&34,43&1,4\%  &&287,43&1,7 \%  \\ 
			\cmidrule(lr){1-1}	\cmidrule(lr){4-4}	\cmidrule(lr){5-5}\cmidrule(lr){7-7}\cmidrule(lr){8-8}	\cmidrule(lr){10-10}\cmidrule(lr){11-11}
			Total/Avg.   &		    && \text{108} &1,1\%& & \text{56} &1,1\% & & \text{386}&1,2\% \\
			\bottomrule		
		\end{tabular}
		\caption{Comparing the plain MILP formulations using the a- and b- instances.}
		\label{tab:comparison_ab}
	\end{threeparttable}
\end{table}

\setlength{\tabcolsep}{3pt}
\sisetup{table-align-text-post = false, table-number-alignment=center,table-format = 4.1}
\begin{table}[ht]
	\centering\footnotesize
	\begin{threeparttable}		
		\begin{tabular}{cSS[table-format = 1.1]SS[table-format = 1.1]SS[table-format = 1.1]SS[table-format = 1.1]SS[table-format = 1.1]SS[table-format = 1.1]}
			\toprule
			& \multicolumn{4}{c}{EB} & \multicolumn{4}{c}{LAEB} & \multicolumn{4}{c}{ALAEB} \\
			\cmidrule(lr){2-5} \cmidrule(lr){6-9} \cmidrule(lr){10-13}
			$\text{Inst.}$ &  $\text{Obj.}$  & $\text{Gap}$  & $\text{CPU}$ & $\text{rGap}$ & $\text{Obj.}$  &  $\text{Gap}$  & $\text{CPU}$ & $\text{rGap}$ &  $\text{Obj.}$ &  $\text{Gap}$  & $\text{CPU}$ & $\text{ rGap}$\\
			\midrule
			a2-16-X &278,2&&0,14&1,3 \%&278,2&&0,08&1,7 \%&278,2&&0,09&1,7 \% \\
			a2-20-X &330,7&&0,03&2,0 \%&330,7&&0,02&2,0 \%&330,7&&0,05&2,4 \% \\
			a2-24-X &389,1&&0,26&2,6 \%&389,1&&0,20&2,6 \%&389,1&&0,62&2,8 \% \\
			a3-18-X &272,7&&0,07&2,5 \%&272,7&&0,11&3,0 \%&272,7&&0,39&7,0 \% \\
			a3-24-X &289,6&&0,73&2,2 \%&289,6&&0,75&3,5 \%&289,6&&1,99&3,5 \% \\
			a3-30-X &452,8&&0,34&5,0 \%&452,8&&0,28&5,4 \%&452,8&&1,32&5,6 \% \\
			a3-36-X &501,0&&0,37&1,7 \%&501,0&&0,29&1,9 \%&501,0&&1,08&1,9 \% \\
			a4-16-X &235,2&&0,43&5,9 \%&235,2&&0,34&5,9 \%&235,2&&1,54&5,9 \% \\
			a4-24-X &359,4&&0,37&2,3 \%&359,4&&0,12&2,5 \%&359,4&&0,71&4,2 \% \\
			a4-32-X &447,3&&6,39&3,8 \%&447,3&&2,54&3,9 \%&447,3&&19,78&4,5 \% \\
			a4-40-X &509,0&&38,56&2,9 \%&509,0&&21,93&3,0 \%&509,0&&85,36&3,2 \% \\
			a4-48-X &620,3&&1194,49&6,4 \%&620,3&&660,95&7,8 \%&622,5&2,5\% &2h&10,6 \%  \\
			a5-40-X &464,0&&16,77&3,1 \%&464,0&&8,79&3,3 \%&464,0&&33,64&4,5 \%  \\
			a5-50-X &621,9&&285,24&5,8 \%&621,9&&213,63&5,5 \%&621,9&0,5\% &2h&6,6 \%  \\
			a5-60-X &745,4&&281,67&3,6 \%&745,4&&254,23&4,0 \%&745,4&0,5\% &2h&6,2 \% \\
			a6-48-X &572,5&&7082,45&5,5 \%&572,5&&2675,00&5,6 \%&572,5&1,1\% &2h&6,4 \% \\
			a6-60-X &757,9&&1993,91&3,7 \%&757,9&&855,54&4,6 \%&757,9&1,8\% &2h&5,5 \% \\
			a6-72-X & N/A & N/A &2h&N/A&868,4&3,4\% &2h&6,8 \%&869,6&4,5\% &2h&8,5 \% \\
			a7-56-X &663,5&&4411,85&5,0 \%&663,5&&5131,52&5,7 \%&668,9&1,9\% &2h&6,7 \% \\
			a7-70-X &815,3&0,6\% &2h&6,2 \%&815,3&&1592,16&6,4 \%&815,3&3,0\% &2h&8,0 \% \\
			a7-84-X & N/A & N/A &2h&N/A& N/A & N/A &2h&N/A& N/A & N/A &2h&N/A \\
			a8-64-X &702,6&1,9\% &2h&5,6 \%&701,2&0,9\% &2h&5,5 \%&701,4&2,8\% &2h&7,1 \% \\
			a8-80-X & N/A & N/A &2h&N/A& N/A & N/A &2h&N/A& N/A & N/A &2h&N/A \\
			a8-96-X & N/A & N/A &2h&N/A& N/A & N/A &2h&N/A& N/A & N/A &2h&N/A \\
			\cmidrule(lr){1-1}	\cmidrule(lr){4-4}\cmidrule(lr){5-5}	\cmidrule(lr){8-8}\cmidrule(lr){9-9}	\cmidrule(lr){12-12}\cmidrule(lr){13-13}
			Total/Avg. &               &          & \text{58514} & 3,9\%&   &    &   \text{47418} &4,3\%  &  &   & \text{86546} &5,4\%\\
			\cdashlinelr{1-13}\addlinespace[2mm]
			b2-16-X & 282,5 && 0,31 & 2,2 \% & 282,5 && 0,40 & 2,4 \% & 282,5 && 0,08 & 2,4 \% \\
			b2-20-X & 323,6 && 0,02 & 3,0 \% & 323,6 && 0,02 & 3,0 \% & 323,6 && 0,04 & 3,0 \%  \\
			b2-24-X & 412,3 && 0,03 & 0,0 \% & 412,3 && 0,02 & 0,0 \% & 412,3 && 0,02 & 0,0 \%  \\
			b3-18-X & 290,4 && 0,06 & 3,0 \% & 290,4 && 0,04 & 3,0 \% & 290,4 && 0,07 & 3,0 \%  \\
			b3-24-X & 363,7 && 0,12 & 0,4 \% & 363,7 && 0,06 & 0,4 \% & 363,7 && 0,28 & 0,5 \%  \\
			b3-30-X & 504,3 && 0,17 & 2,3 \% & 504,3 && 0,20 & 2,3 \% & 504,3 && 0,26 & 2,3 \%  \\
			b3-36-X & 565,9 && 0,10 & 0,9 \% & 565,9 && 0,05 & 0,9 \% & 565,9 && 0,06 & 0,9 \%  \\
			b4-16-X & 289,9 && 0,03 & 2,2 \% & 289,9 && 0,02 & 2,3 \% & 289,9 && 0,02 & 2,3 \%  \\
			b4-24-X & 347,0 && 3,55 & 4,4 \% & 347,0 && 0,95 & 5,0 \% & 347,0 && 2,17 & 5,0 \%  \\
			b4-32-X & 491,0 && 0,07 & 1,0 \% & 491,0 && 0,05 & 1,0 \% & 491,0 && 0,06 & 1,0 \%  \\
			b4-40-X & 628,3 && 0,21 & 1,0 \% & 628,3 && 0,25 & 1,1 \% & 628,3 && 0,37 & 1,1 \%  \\
			b4-48-X & 627,4 && 4,89 & 1,4 \% & 627,4 && 1,61 & 1,4 \% & 627,4 && 6,68 & 2,0 \%  \\
			b5-40-X & 585,1 && 5,53 & 2,7 \% & 585,1 && 1,72 & 2,8 \% & 585,1 && 24,13 & 2,7 \%  \\
			b5-50-X & 708,8 && 2,31 & 0,8 \% & 708,8 && 1,12 & 0,8 \% & 708,8 && 5,33 & 0,8 \%  \\
			b5-60-X & 851,9 && 10,72 & 1,8 \% & 851,9 && 5,03 & 1,9 \% & 851,9 && 25,92 & 2,6 \%  \\
			b6-48-X & 691,6 && 2,58 & 1,7 \% & 691,6 && 0,79 & 1,7 \% & 691,6 && 1,94 & 2,2 \%  \\
			b6-60-X & 841,6 && 2,32 & 1,7 \% & 841,6 && 0,63 & 1,8 \% & 841,6 && 5,14 & 1,8 \%  \\
			b6-72-X & 930,3 && 52,73 & 1,7 \% & 930,3 && 49,56 & 1,5 \% & 930,3 && 31,17 & 1,9 \%  \\
			b7-56-X & 787,9 && 142,77 & 1,4 \% & 787,9 && 4,32 & 1,5 \% & 787,9 && 285,5 & 1,5 \%  \\
			b7-70-X & 865,3 && 19,34 & 2,6 \% & 865,3 && 14,68 & 2,5 \% & 865,3 && 132,21 & 2,8 \%  \\
			b7-84-X & 1141,2 && 3676,19 & 3,2 \% & 1141,2 && 2428,40 & 3,2 \% & 1141,9 & 1,0\%  & 2h & 3,6 \%  \\
			b8-64-X & 818,3 && 12,07 & 1,9 \% & 818,3 && 9,16 & 2,0 \% & 818,3 && 24,68 & 2,4 \%  \\
			b8-80-X & 998,3 && 16,10 & 2,2 \% & 998,3 && 14,04 & 2,5 \% & 998,3 && 27,32 & 2,5 \%  \\
			b8-96-X & 1137,7 & 0,3\%  & 2h & 2,0 \% & 1137,7 && 3425,66 & 2,1 \% & 1139,5 & 1,1\%  & 2h & 2,7 \%  \\
			\cmidrule(lr){1-1}	\cmidrule(lr){4-4}\cmidrule(lr){5-5}	\cmidrule(lr){8-8}\cmidrule(lr){9-9}	\cmidrule(lr){12-12}\cmidrule(lr){13-13}
			Total/Avg.  &			&&\text{11152}& 1,9\%& && \text{5958} &2,0\% & & & \text{14973} &2,1\% \\
			\bottomrule		
		\end{tabular}
		\caption{Comparing the plain MILP formulations using the a-X and b-X instances.}
		\label{tab:comparison_aXbX}
	\end{threeparttable}
\end{table}

\begin{landscape}
	\setlength{\tabcolsep}{3pt}
	\sisetup{table-align-text-post = false, table-number-alignment=center,table-format = 4.1}
	\begin{table}[ht]
		\centering\footnotesize
		\begin{threeparttable}		
			\begin{tabular}{cccccccccSS[table-format = 1.1]cc}
				\toprule
				$\text{GP}$ & $\text{VS1}$ & $\text{VS2}$& $\text{VS3}$ & $\text{VS4}$ & $\text{CI1}$ & $\text{IP1}$ & $\text{IP2}$ & $\text{\# Improvements}$ & $\text{Total CPU}$ & $\text{R. gap}$ & $\text{Avg. Deviation}$ & $\text{Avg.* Deviation}$\\ 
				\midrule
				0&0&0&0&0&0&0&0& $\text{-}$ &407.89& 2.9\% & $\text{-}$ & $\text{-}$ \\
				\midrule
				1&1&1&1&1&1&1&1&12&340,57& 1.9\% &-17 \%&-64 \% \\
				%\cdashlinelr{1-13}
				%0&1&1&1&1&1&0&0&7&453,96& 2.6\%  &11 \%&-57 \% \\
				1&0&1&1&1&1&1&1&12&195,70&1.9\% & -52 \%&-82 \% \\
				1&1&0&1&1&1&1&1&12&196,69&1.9\% & -52 \%&-79 \% \\
				1&1&1&0&1&1&1&1&13&278,88&1.8\% &-32 \%&-68 \% \\
				1&1&1&1&0&1&1&1&13&289,63&1.9\% & -29 \%&-75 \% \\
				1&1&1&1&1&0&1&1&12&306,77&1.9\% & -25 \%&-65 \% \\
				1&1&1&1&1&1&0&1&13&322,33&2.0\% & -21 \%&-70 \% \\
				%1&1&1&1&1&1&1&0&12&197,49&2.2\% & -52 \%&-75 \% \\
				%1&1&1&1&1&1&0&0&11&366,52&2.5\%& -10 \%&-66 \% \\
				\cdashlinelr{1-13}
				1&1&0&0&0&0&0&0&11&299,61&2.6\% & -27 \%&-66 \% \\
				1&0&1&0&0&0&0&0&10&345,21&2.5\% & -15 \%&-46 \% \\
				1&0&0&1&0&0&0&0&9&425,08&2.6\% & 4 \%&-59 \% \\
				1&0&0&0&1&0&0&0&10&390,51&2.7\% & -4 \%&-71 \% \\
				1&0&0&0&0&1&0&0&11&302,21&2.6\% & -26 \%&-60 \% \\
				1&0&0&0&0&0&1&0&12&168,40&2.4\% & -58 \%&-74 \% \\
				1&0&0&0&0&0&0&1&13&307,55&2.1\% & -25 \%&-69 \% \\
				\cdashlinelr{1-13}
				1&0&0&0&1&1&1&1&14&171,31&1.9\% & -58 \%&-75 \%  \\
				\bottomrule		
			\end{tabular}
			\caption{Comparing individual and joint effect of preprocessing components on the test instances. A zero represents a component to be switched off and a one represents it to be switched on.}
			\label{tab:pre}
		\end{threeparttable}
	\end{table}
	\setlength{\tabcolsep}{3pt}
	\sisetup{table-align-text-post = false, table-number-alignment=center,table-format = 4.1}
	\begin{table}[ht]
		\centering\footnotesize
		\begin{threeparttable}		
			\begin{tabular}{cccS[table-format = 1.1]cSS[table-format = 1.1]S[table-format = 1.1]}
				\toprule
				\text{LBsUBs}    & \text{FA}   &\text{FP} & \text{Avg. rGap} &\text{\# Improvements}&\text{Total CPU}&\text{Avg. Deviation}&\text{Avg.* Deviation} \\
				\midrule
				0&0&0&1,9 \%&\text{-}&340,57&\text{-}&\text{-} \\
				\midrule
				1&1&1&1,9 \%&5&583,29&66 \%&-25 \% \\
				0&1&1&1,9 \%&6&557,65&54 \%&-20 \% \\
				1&0&1&1,9 \%&6&512,79&56 \%&-25 \% \\
				1&1&0&1,9 \%&5&580,06&66 \%&-25 \% \\
				1&0&0&1,9 \%&6&508,42&56 \%&-25 \% \\
				0&1&0&1,9 \%&6&555,10&54 \%&-20 \% \\
				\bottomrule		
			\end{tabular}
			\caption{Evaluating the performance of the branch-and-but algorithm. A zero represents a component to be switched off and a one represents it to be switched on.}
			\label{tab:usercuts}
		\end{threeparttable}
	\end{table}
\end{landscape}

\begin{landscape}
	\setlength{\tabcolsep}{5pt}
	\sisetup{table-align-text-post = false, table-number-alignment=center,table-format = 3.1}
	\begin{table}[ht]
		\centering\footnotesize
		\begin{threeparttable}		
			\begin{tabular}{cSS[table-format = 1.2]S[table-format = 1.1] c@{\hspace{6mm}} cSS[table-format = 1.2]S[table-format = 1.1] c@{\hspace{6mm}} cSS[table-format = 1.1]SS[table-format = 1.1]c@{\hspace{6mm}}cSS[table-format = 1.2]S[table-format = 1.1]}
				\toprule
				$\text{Inst.}$ &  $\text{Obj.}$  & $\text{CPU}$ & $\text{rGap}$ && $\text{Inst.}$ &  $\text{Obj.}$ & $\text{CPU}$ & $\text{rGap}$ &&$\text{Inst.}$ &  $\text{Obj.}$  & $\text{Gap}$  & $\text{CPU}$ & $\text{rGap}$ & 	& $\text{Inst.}$ &  $\text{Obj.}$ &  $\text{CPU}$ & $\text{rGap}$ \\
				\midrule
				  a2-16  & 294,3 & 0,10 & 0,0 \% && b2-16 & 309,4 & 0,17 & 0,6 \% &&  a2-16-X  & 278,2 && 0,07 & 0,0 \% && b2-16-X & 282,5 & 0,57 & 1,3 \% \\
				a2-20  & 344,9 & 0,02 & 0,0 \% && b2-20 & 332,7 & 0,02 & 0,0 \% &&  a2-20-X  & 330,7 && 0,03 & 0,0 \% && b2-20-X & 323,6 & 0,02 & 0,8 \% \\ 
				a2-24  & 431,1 & 0,05 & 0,8 \% && b2-24 & 444,7 & 0,07 & 1,0 \% &&  a2-24-X  & 389,1 && 0,13 & 1,4 \% && b2-24-X & 412,3 & 0,03 & 0,0 \%  \\
				a3-18  & 300,5 & 0,02 & 0,3 \% && b3-18 & 301,6 & 0,03 & 0,9 \% &&  a3-18-X  & 272,7 && 0,11 & 2,8 \% && b3-18-X & 290,4 & 0,05 & 1,3 \%  \\
				a3-24  & 344,8 & 0,07 & 1,6 \% && b3-24 & 394,5 & 0,12 & 1,9 \% &&  a3-24-X  & 289,6 && 0,74 & 3,3 \% && b3-24-X & 363,7 & 0,08 & 0,4 \%  \\
				a3-30  & 494,8 & 0,06 & 0,0 \% && b3-30 & 531,4 & 0,03 & 0,0 \% &&  a3-30-X  & 452,8 && 0,27 & 2,2 \% && b3-30-X & 504,3 & 0,16 & 1,6 \%  \\
				a3-36  & 583,2 & 0,17 & 2,5 \% && b3-36 & 603,8 & 0,05 & 0,1 \% &&  a3-36-X  & 501,0 && 0,48 & 0,6 \% && b3-36-X & 565,9 & 0,07 & 0,0 \%  \\
				a4-16  & 282,7 & 0,04 & 0,7 \% && b4-16 & 296,9 & 0,02 & 0,6 \% &&  a4-16-X  & 235,2 && 0,50 & 5,4 \% && b4-16-X & 289,9 & 0,02 & 1,7 \%  \\
				a4-24  & 375,0 & 0,04 & 0,0 \% && b4-24 & 371,4 & 0,04 & 0,1 \% &&  a4-24-X  & 359,4 && 0,15 & 1,5 \% && b4-24-X & 347,0 & 0,47 & 2,1 \%  \\
				a4-32  & 485,5 & 0,12 & 1,1 \% && b4-32 & 494,9 & 0,04 & 0,0 \% &&  a4-32-X  & 447,3 && 1,78 & 3,4 \% && b4-32-X & 491,0 & 0,08 & 0,7 \%  \\
				a4-40  & 557,7 & 0,28 & 1,2 \% && b4-40 & 656,6 & 0,09 & 0,3 \% &&  a4-40-X  & 509,0 && 0,95 & 2,1 \% && b4-40-X & 628,3 & 0,28 & 0,1 \%  \\
				a4-48  & 668,8 & 0,27 & 1,7 \% && b4-48 & 673,8 & 0,65 & 0,7 \% &&  a4-48-X  & 620,3 && 56,36 & 5,5 \% && b4-48-X & 627,4 & 1,46 & 0,9 \%  \\
				a5-40  & 498,4 & 0,16 & 0,6 \% && b5-40 & 613,7 & 0,12 & 0,4 \% & & a5-40-X  & 464,0 && 6,12 & 2,9 \% && b5-40-X & 585,1 & 2,65 & 1,9 \%  \\
				a5-50  & 686,6 & 1,51 & 2,0 \% && b5-50 & 761,4 & 0,21 & 1,0 \% &&  a5-50-X  & 621,9 && 73,02 & 4,0 \% && b5-50-X & 708,8 & 0,53 & 0,6 \%  \\
				a5-60  & 808,3 & 0,45 & 1,4 \% && b5-60 & 902,0 & 0,86 & 1,6 \% & & a5-60-X  & 745,4 && 39,78 & 3,3 \% && b5-60-X & 851,9 & 1,66 & 1,3 \%  \\
				a6-48  & 604,1 & 0,38 & 0,4 \% && b6-48 & 714,8 & 0,18 & 0,3 \% &  &a6-48-X  & 572,5 && 935,75 & 4,7 \% && b6-48-X & 691,6 & 0,53 & 0,4 \%  \\
				a6-60  & 819,3 & 2,14 & 1,6 \% && b6-60 & 860,0 & 0,24 & 0,1 \% &&  a6-60-X  & 757,9 && 332,38 & 3,2 \% && b6-60-X & 841,6 & 1,13 & 1,3 \%  \\
				a6-72  & 916,1 & 11,34 & 2,2 \% && b6-72 & 978,5 & 3,95 & 0,9 \% &&  a6-72-X  & 869,6 & 2,1 \% & 2h & 5,6 \% && b6-72-X & 930,3 & 7,16 & 1,4 \%  \\ 
				a7-56  & 724,0 & 2,32 & 1,2 \% && b7-56 & 824,0 & 1,46 & 0,8 \% & & a7-56-X  & 663,5 && 551,48 & 4,1 \% && b7-56-X & 787,9 & 5,06 & 0,8 \%  \\
				a7-70  & 875,7 & 2,11 & 0,9 \% && b7-70 & 912,6 & 0,65 & 1,5 \% &&  a7-70-X  & 815,3 && 711,65 & 4,7 \% && b7-70-X & 865,3 & 6,11 & 1,6 \%  \\
				a7-84  & 1033,3 & 11,13 & 1,6 \% && b7-84 & 1203,3 & 1,50 & 0,8 \% &&  a7-84-X  & 950,6 & 1,4 \% & 2h & 6,2 \% && b7-84-X & 1141,2 & 89,10 & 2,1 \%  \\
				a8-64  & 747,5 & 2,29 & 1,1 \% &&b8-64 & 839,9 & 0,70 & 1,6 \% &&  a8-64-X  & 701,2 && 3042,70 & 4,1 \% && b8-64-X & 818,3 & 2,38 & 1,8 \%  \\
				a8-80  & 945,8 & 7,58 & 1,9 \% && b8-80 & 1036,4 & 0,64 & 0,8 \% &&  a8-80-X  & 880,3 & 2,7 \% & 2h & 6,4 \% && b8-80-X & 998,3 & 6,00 & 1,7 \% \\ 
				a8-96  & 1229,7 & 425,06 & 2,8 \% && b8-96 & 1185,6 & 11,13 & 1,1 \% &&  a8-96-X  & N/A & N/A & 2h & N/A && b8-96-X & 1137,7 & 340,51 & 2,3 \%  \\
				\cmidrule(lr){1-1} \cmidrule(lr){3-3} \cmidrule(lr){4-4}\cmidrule(lr){8-8}\cmidrule(lr){9-9} \cmidrule(lr){14-14}\cmidrule(lr){15-15} \cmidrule(lr){19-19}\cmidrule(lr){20-20}   
				Total  & & \text{467} &1.2\%&& & & \text{22} &0.7\%&& & & & \text{34554}& 3.4\%& & & & \text{466}&1.2\%\\
				\bottomrule		
			\end{tabular}
			\caption{Results for the benchmark instances solved using the LAEB formulation and preprocessing.}
			\label{tab:1thread_ab_aXbX}
	\end{threeparttable}
	\end{table}
\end{landscape}

\newpage
\section*{Figures}

\newpage
\appendix
\section*{Electronic Companions}

\section{Notation}\label{Notation}
We summarize the notation of this paper in Table~\ref{tab_notation}.
{\small
	\begin{table}[ht]
		\begin{tabular}{ll}
			%$A$ & set of arcs in the event-based graph \\
			$B_v$ & variable in the EB formulation indicating start of service at location $v_1$\\
					& if event $v$ occurs \\
					$B_v^{LB}$ & lower bound for $B_v$ \\
					$B_v^{UB}$ & upper bound for $B_v$ \\
					$B_{vw}^{LB}$ & earliest feasible beginning of service at node $w$, when $x_{(v,w)}=1$\\
					$B_{vw}^{UB}$ & latest feasible beginning of service at node $v$, when $x_{(v,w)}=1$\\
					$\bar{B}_i^{UB}$ & $= \max_{v\in V_{i^+}} \{B_v^{UB}\}$, i.e., upper bound on beginning of service at location $i^+$\\
					$\bar{B}_j$ & variable indicating start of service at location $j\in J$ or $j \in \bar{J}$, respectively \\
					$\bar{c}_{ij}$ & costs for travelling between locations $i, j \in J$ or $i, j \in \bar{J}$, respectively \\
					$c_a$ & costs for travelling arc $a = (v,w)\in A$, i.e., between locations $v_1$ and $w_1$ \\
					$D$ & set of delivery locations, i.e., $D = \{1^-, \ldots, n^- \}$ \\
					$\delta^{in}(v)$ & $= \{ (w,v) \in A \}$, i.e., ingoing arcs of node $v$ \\
					$\delta^{out}(v)$ & $= \{ (v,w) \in A \}$, i.e., outgoing arcs of node $v$ \\
					$e_j$ & earliest possible start of service at location $j\in J$ or $j \in \bar{J}$, respectively \\
					$f^1_{ij}$, $f^2_{ij}$ & feasibility of the paths $j^+ \rightarrow i^+ \rightarrow j^- \rightarrow i^-$ and $j^+ \rightarrow i^+ \rightarrow i^- \rightarrow j^-$,\\
					& respectively, w.r.t.\ ride time and time window constraints, $i,j\in R$ \\
					$J$ & set of all locations in the event-based and location-augmented-event-based \\
					& formulations, i.e., $J=P\cup D \cup\{ 0\}$ \\ %\ $J = \{0, 1^+, \ldots, n^+, 1^-, \ldots, n^- \}$ \\
					$\bar{J}$ & set of locations in the LB formulation, \\
					& i.e., $\bar{J}= J\setminus \{0\} \cup{0^+,0^-}$\\
					$i^+$, $i^-$& pick-up and delivery location of customer $i$ \\
					$K$ & number of vehicles \\
%					$u^m_{i^+}$ & unique event s.t.\ $(u^m_{i^+})_1 = i^+$ and all users sitting in the vehicle directly \\
%					& after event $u^m$ are in $u^m_{i^+}$ in the vehicle too \\
%					$u^m_{i^- }$ & unique event s.t.\ $(u^m_{i^-})_1 = i^-$ and all users sitting in the vehicle directly \\
%					& before event $u^m$ are in $u^m_{i^-}$ in the vehicle too \\
					$L_i$ & maximum ride time of customer $i$ \\
			$\ell_j$ & latest possible start of service at location $j\in J$ or $j \in \bar{J}$, respectively \\
					$\bar{M}_{ij}, M_{vw}$ & sufficiently large constants \\
					$n$ & number of customers \\
					$P$ & set of pick-up locations, i.e., $P = \{ 1^+, \ldots, n^+ \}$ \\
					$Q$ & vehicle capacity \\
					$Q_j$ & variable in the LB model indicating the vehicle load after \\
					& visiting location $j\in\bar{J}$ \\
					%$Q_j^{(u,v)}$ & $= x_{(v,w)} \sum_{j=1}^Q q_j \cdot w_j$ if $j \in P$ and $= x_{(v,w)} \sum_{j=2}^Q q_j \cdot w_j$ if $j \in D$ \\
					$q_j$ & number of persons entering or leaving the vehicle at location $j\in J$ \\
					& or $j\in \bar{J}$, respectively \\
					$R$ & set of all requests, i.e., $R = \{ 1, \ldots, n \}$ \\
					$S \in \mathcal{S}$ & set for pairing and precedence constraints \eqref{eq:SEC} \\
					%$\bar{S}$ & sequence of events $\bar{S} = \{(u^1,u^2), \ldots, (u^{m-1},u^m)\}$, $m \in \mathbb{N}$ arbitrary but fixed  \\
					$s_j$ & service time at location $j\in J$ or $j\in \bar{J}$, respectively \\
					$T$ & service duration \\
					$TW$ & $= \ell_j - e_j$ length of pick-up (inbound requests) or delivery (outbound \\ & requests) time window \\
					$\bar{t}_{ij}$ & travel time between locations $i, j \in J$ or $i, j \in \bar{J}$, respectively \\
					$\bar{t}_i$ & travel time for a direct ride from pick-up to delivery location of request $i$ \\
					$t_{(v,w)}$ & travel time between locations $v_1$ and $w_1$ \\
					%$u,v,w,w' \in V$ & events \\
					$V$ & set of nodes in event-based graph \\
					$V_0$ & event $(0, \ldots, 0)$ \\
					$V_{i^+}$ & set of nodes $v\in V$ with $v_1 = i^+$ \\
					$V_{i^-}$ & set of nodes $v\in V$ with $v_1 = i^-$ \\
					%$(v_1, \ldots, v_Q)$ & occupancy of the in event $v$ where $v_1$ is the current location and $v_2, \ldots, v_Q$  \\
					& are in decreasing order (remaining positions are 0 if less than $Q$ customers \\
					& sitting in the vehicle) \\
					$x_a$ & binary variable in EB and LAEB model, continuous variable in ALAEB model \\
					$\bar{x}_{ij}$ & binary variable in LB and ALAEB model\\
%					$Z_{(\cdot)}$ & objective value \\
%					$Z_{(\cdot)}^{rel}$ & objective value of the LP relaxation 
		\end{tabular}		
		\caption{Notation used in the paper} 
		\label{tab_notation}
\end{table}}

\section{Proof of Theorem \ref{equi}} \label{PTequi}
Constraints \eqref{eq:B}, \eqref{eq:mrt}, and \eqref{E:flow}--\eqref{E:K} are part of both model formulations. As we consider the LP relaxation, constraints \eqref{E:binary} and \eqref{barxrel} are identical. Thus, we have to show that \eqref{E2:objective}--\eqref{E2:time} are equivalent to \eqref{eq:xbin}, \eqref{xbarx}, and \eqref{E3:objective}--\eqref{EL2:time}.

Starting with \eqref{eq:xbin}, \eqref{xbarx}, and \eqref{E3:objective}--\eqref{EL2:time}, we can replace $\bar{x}_{ij}$ in objective \eqref{E3:objective} and constraints \eqref{EL2:time} by $\sum_{(v,w) \in A: v_1 = i \wedge w_1 = j} x_{(v,w)}$ (equations \eqref{xbarx}) and obtain objective \eqref{E2:objective} and constraints \eqref{E2:time}, respectively. Moreover, equations \eqref{xbarx} and constraints \eqref{E:flow} and \eqref{E:pick} imply $0 \leq \bar{x}_{ij} \leq 1$ for all $i,j \in \bar{J}$. Thus, constraints \eqref{eq:xbin} are redundant in the LP relaxation of ALAEB. However, then $\bar{x}_{ij}$ is only set in \eqref{xbarx} but not used in the model any more such that \eqref{xbarx} is also redundant and we obtain \eqref{E2:objective}--\eqref{E2:time}.

Analogously, we can start with \eqref{E2:objective}--\eqref{E2:time}, add the redundant equations \eqref{xbarx} and constraints \eqref{eq:xbin}. Using equations \eqref{xbarx}, we can replace $\sum_{(v,w) \in A: v_1 = i \wedge w_1 = j} x_{(v,w)}$ by $\bar{x}_{ij}$ in objective \eqref{E2:objective} and constraints \eqref{E2:time} to obtain objective \eqref{E3:objective} and constraints \eqref{EL2:time}, respectively. By this, we obtain \eqref{eq:xbin}, \eqref{xbarx}, and \eqref{E3:objective}--\eqref{EL2:time}.

\section{Proof of Theorem \ref{Net}}\label{ProofNet}
%\begin{Proof}

Condition \ref{C1} ensures that maximum ride time constraints \eqref{eq:mrt} are fulfilled. Because of the first case in condition \ref{C2}
\[
\bar{B}_i + s_i + \bar{t}_{ij} \leq l_i + s_i + \bar{t}_{ij} \leq e_j \leq \bar{B}_j,
\]
i.e., $\bar{B}_i + s_i + t_{ij} \leq \bar{B}_j$ holds for these pairs $i,j$. In the second case, $e_i + s_i + \bar{t}_{ij} > l_j$, location $j$ is not reached within the time window even if the vehicle starts as early as possible in location $i$. Thus, $x_{(v,w)} = 0$ for all $(v,w)$ with $v_1 = i$ and $w_1 = j$. In this case, we can simply omit arc $(v,w)$ from the graph. Due to condition \ref{C2} all remaining arcs fulfill $l_i + s_i + \bar{t}_{ij} \leq e_j$, i.e., location $j$ is reached latest at the beginning of the time window if the vehicle visits location $i$ directly before. Thus, $\bar{B}_i + s_i + t_{ij} \leq \bar{B}_j$ with $v_1 = i$ and $w_1 = j$ is always fulfilled. As a consequence constraints \eqref{E2:time} are always fulfilled. For the same reason we always find a solution with $\bar{B}_i \in [e_i,l_i]$ for all $i \in \bar{J}$. Hence, we can always find a solution fulfilling the time window constraints \eqref{eq:B}. Together, we can omit constraints \eqref{eq:B}, \eqref{eq:mrt}, and \eqref{E2:time}.

%If $e_{j}= \ell_j$, $\bar{B}_j= e_j$ due to constraints \eqref{eq:B}. Then, maximum ride time constraints \eqref{eq:mrt} are fulfilled or not, i.e., we directly know whether the solution is feasible due to the maximum ride time. To ensure that \eqref{E2:time} is fulfilled, we delete all arcs $(v,w)$ with $v_1 = i$ and $w_1 = j$ if $\bar{B}_{j} < \bar{B}_{i} + s_{i} + \bar{t}_{ij}$. Thus, \eqref{E2:time} is fulfilled if $x_{(v,w)} = 1$. Due to the choice of big M, it is always fulfilled if $x_{(v,w)} = 0$. Together, we can omit constraints \eqref{eq:B}, \eqref{eq:mrt}, and \eqref{E2:time}. 

Constraints \eqref{E:K} are equivalent to $\sum_{a \in \delta^{out}(0, \ldots, 0)} x_a = K$ by $c_{(0, \ldots, 0),(0, \ldots, 0)} = 0$, i.e., vehicles can drive from the depot to the depot without any costs. Without constraints \eqref{E:pick}, the remaining formulation \eqref{E:flow}, \eqref{E:K}, and \eqref{E:binary} with objective \eqref{E2:objective} is a Minimum Cost Flow Problem with arc capacities $x_a \leq 1$.
\begin{figure}[htbp]
	\centering
	\scalebox{0.8}{\begin{tikzpicture}
\node at (9.5,5) {$(w^m,v^m) \in \delta^{in}(v^m), v^m_1 = i^+ \;\forall m = 1, ..., \bar{m}$};
\node (D) at (5,0) {$w^{\bar{m}}$};
\node (C) at (5,1) {$\vdots$};
\node (B) at (5,2) {$w^2$};
\node (A) at (5,3) {$w^1$};

\node (E) at (8,1.5) {$in_i$};
\node (F) at (11,1.5) {$out_i$};

\node (G) at (14,3) {$v^1$};
\node (H) at (14,2) {$v^2$};
\node (I) at (14,1) {$\vdots$};
\node (J) at (14,0) {$v^{\bar{m}}$};

% \draw[->] (A) to node[above] {\;\;\;\;\;\;\;\;\;\;\;\;\;\;\;\;\;\;$x_{(w^1,in_i)} = x_{(w^1,v^1)}$} (E);
\draw[->] (A) to node[above] {\;\;\;\;$x_{(w^1,in_i)}$} (E);
\draw[->] (B) -- (E);
\draw[->,dotted] (C) -- (E);
\draw[->] (D) -- (E);
\draw[->] (E) to node[above]{$x_{(in_i,out_i)} = 1$} (F);
% \draw[->] (F) to node[above] {$x_{(out_i,v^1)} = x_{(w^1,v^1)}$\;\;\;\;\;\;\;\;\;\;\;\;\;\;\;\;\;\;\;\;\;\;} (G);
\draw[->] (F) to node[above] {$x_{(out_i,v^1)}$\;\;\;\;\;\;\;} (G);
\draw[->] (F) -- (H);
\draw[->,dotted] (F) -- (I);
\draw[->] (F) -- (J);

\node (D1) at (0,0) {$w^{\bar{m}}$};
\node (C1) at (0,1) {$\vdots$};
\node (B1) at (0,2) {$w^2$};
\node (A1) at (0,3) {$w^1$};
\node (G1) at (3,3) {$v^1$};
\node (H1) at (3,2) {$v^2$};
\node (I1) at (3,1) {$\vdots$};
\node (J1) at (3,0) {$v^{\bar{m}}$};
\draw[->] (A1) to node[above]{$x_{(w^1,v^1)}$} (G1);
\draw[->] (B1) to node[above]{$x_{(w^2,v^2)}$} (H1);
\draw[->,dotted] (C1) to (I1);
\draw[->] (D1) to node[above]{$x_{(w^{\bar{m}},v^{\bar{m}})}$} (J1);

\node (D2) at (16,0) {$w^{\bar{m}}$};
\node (C2) at (16,1) {$\vdots$};
\node (B2) at (16,2) {$w^2$};
\node (A2) at (16,3) {$w^1$};
\node (G2) at (19,3) {$v^1$};
\node (H2) at (19,2) {$v^2$};
\node (I2) at (19,1) {$\vdots$};
\node (J2) at (19,0) {$v^{\bar{m}}$};
\draw[->] (A2) to node[above]{$x_{(w^1,v^1)}$} (G2);
\draw[->] (B2) to node[above]{$x_{(w^2,v^2)}$} (H2);
\draw[->,dotted] (C2) to (I2);
\draw[->] (D2) to node[above]{$x_{(w^{\bar{m}},v^{\bar{m}})}$} (J2);
\node at (17.5,4) {$\sum_{m=1}^{\bar{m}} x_{(w^m,v^m)} = 1$};
\end{tikzpicture}}
	\caption{Network transformation}
	\label{fig:PNF}
\end{figure}
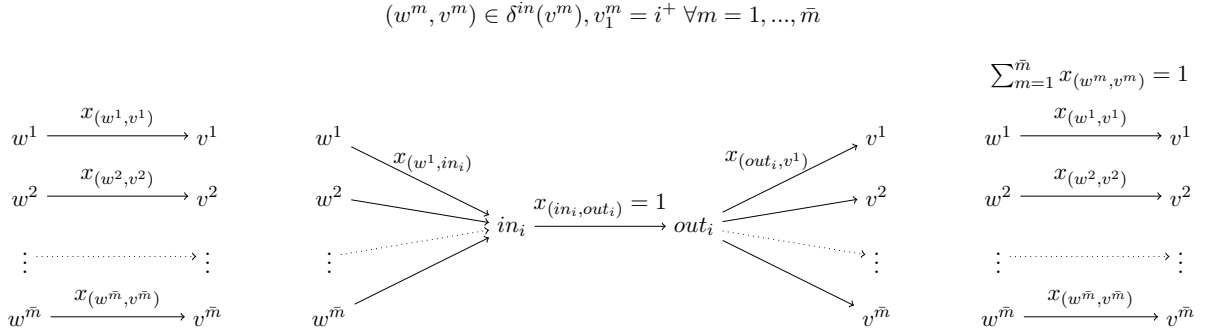

The left side of Figure \ref{fig:PNF} corresponds to \eqref{E:flow}, \eqref{E:K}, and \eqref{E:binary} with objective \eqref{E2:objective} for a fixed location $i^+$, where $v^{m}$, $m = 1, \ldots, \bar{m}$, are all events with $v^{m}_1 = i^+$ and $w^{m}$ the possible predecessor events. Note that $v^m = v^{m'}$ for $m' \neq m$ might hold such that all relations between predecessors $w$ and events $v$ with $v_1 = i^+$ are included. We can add two further nodes $in_i$ and $out_i$ in between and define $out_i$ as a source with an outflow of 1 and $in_i$ as a sink with inflow $1$ without destroying the network flow property. In fact, this is equivalent to set the flow variable $x_{(in_i,out_i)} = 1$ (middle part of the figure). 

As $x_{(w^m,v^m)} = x_{(w^m,in_i)} = x_{(out_i,v^m)}$, 
there is a flow of $x_{(w^m,v^m)}$ from $w^m$ via $in_i$ and $out_i$ to $v^m$. First, this means that $\sum_{m = 1}^{\bar{m}} x_{(w^m,v^m)} = x_{(in_i,out_i)} = 1$. Second, we can shrink $x_{(w^m,in_i)}$, $x_{(out_i,v^m)}$, and the corresponding flow between $in_i$ and $out_i$ to variable $x_{(w^m,v^m)}$ and obtain the situation on the right side of Figure \ref{fig:PNF}.
%\todo{zum Verständnis: diese Lsg ist immer noch Lsg eines MCF Problems weil wir die MCF Eigenschaft durch Hinzufügen von Quelle $out_i$ und Senke $in_i$ mit Bedarf 1 nicht verändert haben. Dadurch ergibt sich Situation auf rechter Seite der Grafik. Dieses Optimierungsproblem entspricht aber LAEB model unter den Bedingungen des Theorems} 
By definition $\sum_{m = 1}^{\bar{m}} x_{(w^m,v^m)}$ sums up all incoming flow to all events with $v_1 = i^+$. Thus, it is equivalent to \eqref{E:pick} for the considered $i \in R$. By repeating this step for all $i \in R$, formulation \eqref{eq:B}, \eqref{eq:mrt}, \eqref{E:flow}--\eqref{E:K}, and \eqref{E:binary}--\eqref{E2:time} is integral if conditions \ref{C1} and \ref{C2} are fulfilled and all arcs $(v,w)$ with $v_1 = i$ and $w_1 = j$ are deleted if $l_j < e_i + s_i + \bar{t}_{ij}$.
%\end{Proof}

\section{Example Location-Based Model}\label{ExlocMod}
Theorem \ref{Net} cannot be transferred to the LB model formulation \eqref{eq:objFunc}--\eqref{eq:xbin} as Example \ref{locNoNet} shows.

\begin{Example} \label{locNoNet} 
	Consider the instance with three customers, $q_1 = q_2 = 2$, $q_3 = 3$, $Q = 6$, and the travel times given in the left part of Figure \ref{fig:Ex_loc}. We assume $\bar{c}_{ij} = \bar{t}_{ij}$ for all $i,j \in \bar{J}$ and that maximum ride time constraints are fulfilled. Moreover, let $e_{1^+} = \ell_{1^+} = 10$, $e_{2^+} = \ell_{2^+} = 11$, $e_{3^+} = \ell_{3^+} = 12$, $e_{1^-} = \ell_{1^-} = 13$, $e_{2^-} = \ell_{2^-} = 14$, $e_{3^-} = \ell_{3^-} = 15$. 
	
	\begin{figure}[htbp]
		\centering
		\scalebox{0.6}{\begin{tikzpicture}
\node at (2,14) {travel times};
\node (D6) at (0,6) {$3^-$};
\node (D5) at (0,9) {$2^-$};
\node (D4) at (0,12) {$1^-$};
\node (D3) at (4,12) {$3^+$};
\node (D2) at (4,9) {$2^+$};
\node (D1) at (4,6) {$1^+$};
\node (D0+) at (2.5,0) {$0^+$};
\node (D0-) at (1.5,0) {$0^-$};

\draw[->] (D0+) to node[above]{10\;\;\;} (D1);
\draw[->] (D0+) to[bend right = 30] node[below]{\;\;\;11} (D2);
\draw[->] (D0+) to[bend right = 40] node[below]{\;\;\;\;12} (D3);
\draw[<->] (D1) to node[above]{1\;\;\;} (D2);
\draw[<->] (D2) to node[above]{1\;\;\;} (D3);

\draw[->] (D6) to node[above]{\;\;\;10} (D0-);
\draw[->] (D5) to[bend right = 30] node[below]{11\;\;\;} (D0-);
\draw[->] (D4) to[bend right = 40] node[below]{12\;\;\;\;} (D0-);
\draw[<->] (D4) to node[above]{\;\;\;1} (D5);
\draw[<->] (D5) to node[above]{\;\;\;1} (D6);

\draw[<->] (D1) to node[above]{1} (D6);
\draw[<->] (D1) to node[below]{2} (D5);
\draw[->] (D1) to node[below]{3} (D4);
\draw[<->] (D2) to node[above]{2} (D6);
\draw[->] (D2) to node[above]{\;\;\;\;\;\;\;\;\;\;\;\;\;\;1} (D5);
\draw[<->] (D2) to node[below]{2} (D4);
\draw[<->] (D3) to node[above]{3} (D6);
\draw[<->] (D3) to node[above]{2} (D5);
\draw[->] (D3) to node[above]{1} (D4);

\draw[<->] (D1) to[bend right = 30] node[above]{2\;\;\;} (D3);
\draw[<->] (D4) to[bend right = 30] node[above]{\;\;\;2} (D6);

\node at (12,14) {location-based model};
\node at (12,13.5) {$Z_L = 24.916$};
\node (L6) at (10,6) {$3^-$};
\node (L5) at (10,9) {$2^-$};
\node (L4) at (10,12) {$1^-$};
\node (L3) at (14,12) {$3^+$};
\node (L2) at (14,9) {$2^+$};
\node (L1) at (14,6) {$1^+$};
\node (L0+) at (12.5,0) {$0^+$};
\node (L0-) at (11.5,0) {$0^-$};

\draw[->] (L0+) to node[above]{1\;\;\;} (L1);
\draw[->] (L0+) to[bend right = 30] node[below]{\;\;\;\;\;\;\;\;\;0.083} (L2);
\draw[->] (L1) to node[above]{0.917\;\;\;\;\;\;\;\;\;\;} (L2);
\draw[->] (L1) to[bend right = 30] node[below]{\;\;\;\;\;\;\;\;\;0.083} (L3);
\draw[->] (L2) to node[above]{0.917\;\;\;\;\;\;\;\;\;\;} (L3);
\draw[->] (L2) to node[above]{0.083} (L5);
\draw[->] (L3) to node[above]{1} (L4);
\draw[->] (L4) to node[above]{\;\;\;\;\;\;\;\;\;\;0.917} (L5);
\draw[->] (L4) to[bend right = 30] node[below]{0.083\;\;\;\;\;\;\;\;\;} (L6);
\draw[->] (L5) to node[above]{\;\;\;\;\;\;\;\;\;\;0.917} (L6);
\draw[->] (L5) to[bend right = 30] node[below]{0.083\;\;\;\;\;\;\;\;\;} (L0-);
\draw[->] (L6) to node[above]{\;\;\;1} (L0-);

\node at (22,14) {location-based model};
\node at (22,13.5) {optimal integer solution};
\node at (22,13) {$Z_L = 46$};
\node (C6) at (20,6) {$3^-$};
\node (C5) at (20,9) {$2^-$};
\node (C4) at (20,12) {$1^-$};
\node (C3) at (24,12) {$3^+$};
\node (C2) at (24,9) {$2^+$};
\node (C1) at (24,6) {$1^+$};
\node (C0+) at (22.5,0) {$0^+$};
\node (C0-) at (21.5,0) {$0^-$};

\draw[->] (C0+) to node[above]{1\;\;\;} (C1);
\draw[->] (C0+) to[bend right = 30] node[below]{\;\;\;1} (C2);
\draw[->] (C1) to[bend right = 30] node[below]{\;\;\;1} (C3);
%\draw[->] (C1) to node[below]{0.01\;\;\;\;\;} (C4);
%\draw[->] (C2) to node[above]{0.01\;\;\;\;\;\;\;} (C3);
\draw[->] (C2) to node[above] {1} (C5);
\draw[->] (C3) to node[above] {1} (C4);
%\draw[->] (C4) to node[below]{\;\;\;\;\;\;\;\;0.01} (C5);
\draw[->] (C4) to[bend right = 30] node[below]{1\;\;\;} (C6);
%\draw[->] (C4) to[bend right = 60] node[below]{0.01\;\;\;\;\;\;\;} (C0);
%\draw[->] (C5) to node[below]{\;\;\;\;\;\;\;\;0.01} (C6);
\draw[->] (C5) to[bend right = 30] node[below]{1\;\;\;} (C0-);
\draw[->] (C6) to node[above]{\;\;\;1} (C0-);

\end{tikzpicture}}
		\caption{Example for LB formulation}
		\label{fig:Ex_loc}
	\end{figure}
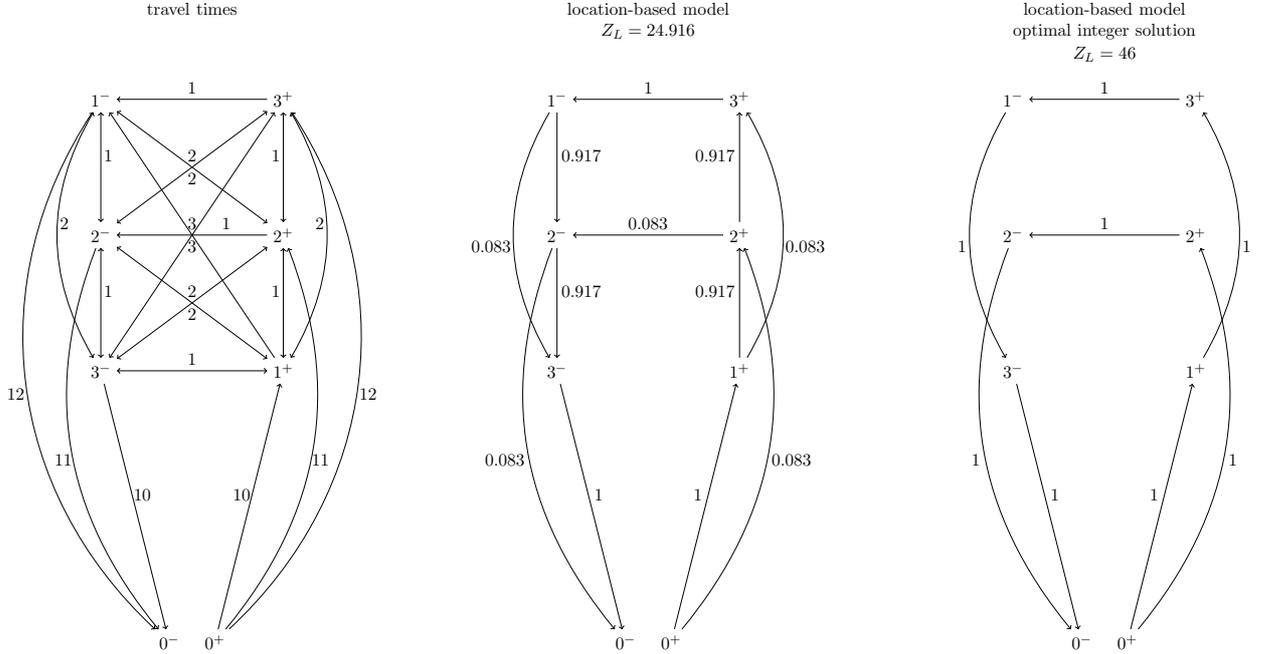
	
	Figure \ref{fig:Ex_loc} presents an optimal solution for \eqref{eq:objFunc}--\eqref{eq:mrt} and $0 \leq \bar{x}_{ij} \leq 1$ for all $i,j \in \bar{J}$ in the middle of the figure and an optimal solution for \eqref{eq:objFunc}--\eqref{eq:xbin} on the right side of the figure. As can be seen there is a positive flow of 0.917 from $1^+$ via $2^+$ to $3^+$, which means that all three customers are in the vehicle at the same time at least for a fractional flow. As they require seven seats in total, but the vehicle has only six, all three customers cannot share the vehicle at the same time. 
\end{Example}

\section{Proof of Theorem \ref{T1}}\label{PT1}
%\begin{Proof}
The proof consists of three steps. First, we prove that every feasible LP solution of the LAEB formulation is also a feasible LP solution of the LB formulation, i.e., $Z_{LAEB}^{rel} \geq Z_{LB}^{rel}$. Second, we give an example that there is a feasible LP solution for the LB formulation which is not LP feasible for the LAEB formulation. Third, we give a concrete instance in which $Z_{LAEB}^{rel} > Z_{LB}^{rel}$ holds.
\begin{enumerate}
	\item Let $\textbf{x} = (x_1, \ldots,x_{|A|})$ be a feasible LP solution for the LAEB formulation. Set $\bar{x}_{ij} = \sum_{v:v_1=i,w:w_1=j} x_{(v,w)}$ for all $i, j \in J$. Furthermore, define $\bar{x}_{0^+j} = \bar{x}_{0^-j} = \bar{x}_{0j}$ and $\bar{x}_{i0^+}= \bar{x}_{i0^-} = \bar{x}_{i0}$. Due to \eqref{E:pick} the ingoing flow is 1 for each pick-up location $i^+$. Because of \eqref{E:flow} this holds also for the outgoing flow. All outgoing arcs of a node $v$ with $v_1=i^+$ end in a node $w$ with $w_j=i$ for one $j = 2, \ldots, Q$ until a node $w$ with $w_1 = i^-$ is reached. This holds also for all nodes connected by an arc with node $w$. With the same argument this holds for their connected nodes and so on.
	%until a node with $v_1=i^-$ is reached. 
	On the other hand there is no arc between a node $w'$ with $w'_1 \neq i^+$  and $w'_l \neq i$ for all $l = 2, \ldots, Q$ and a node $w$ with $w_1=i^-$. Thus, constraints \eqref{E:flow} lead to
	\begin{equation} \label{eq1}
		\sum_{(v,w):w_1=i^+} x_{(v,w)} = 1 = \sum_{(v,w):w_1=i^-} x_{(v,w)}. 
	\end{equation}
	Together, \eqref{eq:degIn} and \eqref{eq:degOut} are fulfilled. With our choice of $\bar{x}_{ij}$, \eqref{E:K} leads to
	\[
	\sum_{j \in P} \bar{x}_{0^+j} = \sum_{(v,w): v_1 = 0} x_{(v,w)}\leq K,
	\]
	which means that \eqref{eq:depot} is fulfilled.
	
	There is one constraint in \eqref{eq:SEC} for each $S \in \mathcal{S}$. Each $S$ consists of a path starting in the starting depot $0^+$ and visiting the delivery location of a request $i$ without having visited the respective pick-up location before. In the EB formulation, all flow starts in the depot \eqref{E:K}, as for all other nodes the flow constraint \eqref{E:flow} holds. Besides, a flow of 1 enters an event-node associated with the pick-up location of $i$ due to \eqref{E:pick}. Thus, and due to our definition of $\bar{x}_{ij}$, a flow of 1 has to leave set $S$. Moreover, due to construction of the event-based graph and our definition of $\bar{x}_{ij}$, all flow leaving a pick-up location of $i$ has to go to a delivery location of $i$, i.e., enter set $S$ again. In total, the flow amongst nodes in $S$ cannot exceed $|S|-2$. Thus, constraints \eqref{eq:SEC} are fulfilled. 
	
	Furthermore, constraints \eqref{eq:cap} and \eqref{eq:q} are fulfilled by $Q_{0^+} = Q_{0^-} = 0$ and
	\[
	Q_j= \id_{\{ j \in P \}}q_j + \sum_{(u,v) \in A: v_1 = j} x_{(u,v)} \sum_{l=2}^Q q_{v_l} \hspace{2cm} \forall j \in P \cup D
	\]
	whereat $\id_{\{ j \in P \}}$ is the indicator function, i.e.,
	\[
	\id_{\{ j \in P \}} = \begin{cases} 1 & \text{if } j \in P \\ 0 & \text{otherwise.} \end{cases}
	\]
	As all terms are non-negative, $Q_j \geq 0$ for all $j \in \bar{J}$. For pick-up locations $Q_j\geq q_j$; thus, $Q_j \geq \max \{ 0,\mathcal{I}_{\{ j \in P \}}q_j \}$ for all $j \in \bar{J}$. Due to construction of event nodes,
	\[
	\sum_{l=1}^Q q_{v_l} \leq Q \hspace{2cm} \forall v \in V
	\]
	and
	\[
	\sum_{l=2}^Q q_{v_l} \leq Q + \mathcal{I}_{\{ j \in P \}}q_j \hspace{2cm} \forall v \in V
	\]
	for delivery locations. Thus, $Q_j \leq \min\{Q,Q+ \mathcal{I}_{\{ j \in P \}}q_j \}$ for all $j \in \bar{J}$ and constraints \eqref{eq:q} are fulfilled. 
	
	Due to construction of the event-based graph, there are no arcs $(v,w)$ with $v_1 = i$ and $w_1 = j$ and
	\[
	\id_{\{ i \in P \}}q_{i} + \sum_{l=2}^Q q_{v_l} + \mathcal{I}_{\{ j \in P \}}q_j > \id_{\{ j\in P \}}q_{j} + \sum_{l=2}^Q q_{w_l}
	\]
	for any pair $i,j \in \bar{J}$. Thus
	\[
	\id_{\{ i\in P \}}q_{i} + \sum_{l=2}^Q q_{v_l} + \mathcal{I}_{\{ j \in P \}}q_j \leq \id_{\{ {j} \in P \}}q_{j} + \sum_{l=2}^Q q_{w_l}
	\]
	holds for every $i,j \in \bar{J}$ if $x_{(v,w)} > 0$ with $v_1 = i$ and $w_1 = j$ and thereby
	\begin{align} \nonumber
		&\sum_{(v,w) \in A: v_1 = i \wedge w_1 = j} x_{(v,w)} \left( \id_{\{ i \in P \}}q_{i} + \sum_{l=2}^Q q_{v_l} + \mathcal{I}_{\{ j \in P \}}q_j \right) \\ \label{xpos}
		\leq &\sum_{(v,w) \in A: v_1 = i \wedge w_1 = j} x_{(v,w)} \left( \id_{\{ j \in P \}}q_{j} + \sum_{l=2}^Q q_{w_l} \right)
	\end{align}
	for every $i, j \in \bar{J}$. Moreover,
	\[
	\id_{\{ i \in P \}}q_{i} + \sum_{l=2}^Q q_{v_l} + \mathcal{I}_{\{ j \in P \}}q_j - Q \leq \id_{\{ j \in P\}} q_{j} \hspace{1cm} \forall i,j \in \bar{J}, v \in V: v_1 = i
	\]
	holds because of construction of the events. Thus,
	\begin{align} \nonumber
		&\Biggl( 1 - \sum_{(v,w) \in A: v_1 = i \wedge w_1 = j} x_{(v,w)} \Biggr) \Biggl( \id_{\{ i \in P \}}q_{i} + \sum_{l=2}^Q q_{v_l} + \mathcal{I}_{\{ j \in P \}}q_j - Q \Biggr) \\ \label{x0}
		\leq &\Biggl( 1 - \sum_{(v,w) \in A: v_1 = i \wedge w_1 = j} x_{(v,w)} \Biggr) \id_{\{ j \in P\}} q_{j} \;\;\; \forall i,j \in \bar{J} .
	\end{align}
	Together, we get
	\[
	\begin{aligned}
		&Q_{i}+ \mathcal{I}_{\{ j \in P \}}q_j - Q \Biggl( 1 - \sum_{(v,w) \in A: v_1 = i \wedge w_1 = j} x_{(v,w)} \Biggr) \\
		\stackrel{\phantom{\eqref{xpos}-\eqref{x0}}}{=} &\id_{\{ i \in P \}}q_{i} + \underbrace{\sum_{(u,v) \in A: v_1 = i} x_{(u,v)}}_{\leq 1 \; \eqref{eq1}} \sum_{l=2}^Q q_{v_l} + \mathcal{I}_{\{ j \in P \}}q_j - Q\Biggl( 1 - \sum_{(v,w) \in A: v_1 = i \wedge w_1 = j} x_{(v,w)} \Biggr) \\
		\stackrel{\phantom{\eqref{xpos}-\eqref{x0}}}{\leq} &\id_{\{ i \in P \}}q_{i} + \sum_{l=2}^Q q_{v_l} + \mathcal{I}_{\{ j \in P \}}q_j - Q\Biggl( 1 - \sum_{(v,w) \in A: v_1 = i \wedge w_1 = j} x_{(v,w)} \Biggr) \\
		\stackrel{\eqref{xpos}-\eqref{x0}}{\leq} &\Biggl( 1 - \sum_{(v,w) \in A: v_1 = i \wedge w_1 = j} x_{(v,w)} \Biggr) \id_{\{ j \in P\}} q_{j} + \sum_{(v,w) \in A: v_1 = i \wedge w_1 = j} x_{(v,w)} \Biggl( \id_{\{ j \in P \}}q_{j} + \sum_{l=2}^Q q_{w_l} \Biggr) \\
		\stackrel{\phantom{\eqref{xpos}-\eqref{x0}}}{\leq} &\id_{\{ j \in P \}}q_{j} + \sum_{(u,w) \in A: w_1 = j} x_{(u,w)} \sum_{l=2}^Q q_{w_l} \\
		\stackrel{\phantom{\eqref{xpos}-\eqref{x0}}}{=} &Q_{j} \hspace{0.5cm} \forall i,j \in \bar{J}.
	\end{aligned}
	\]
	With our choice of $\bar{x}_{ij}$ constraints \eqref{eq:cap} are fulfilled. Together, the first step follows.
	\item
	\begin{figure}[htbp]
		\centering
		\scalebox{0.9}{\begin{tikzpicture}
%bottom left
\node (A) at (-1,0) {$(1^+,0, ..., 0)$};
\node (B) at (3,0) {$(1^-,0, ..., 0)$};
\node (C) at (3,2) {$(2^+,0, ..., 0)$};
\node (D) at (1,3.5) {$(2^-,0, ..., 0)$};
\node (E) at (-1,2) {$(0, ..., 0)$};
\draw[->] (A) to node[below] {1} (B);
\draw[->] (B) to node[right] {1} (C);
\draw[->] (C) to node[above] {1} (D);
\draw[->] (D) to node[above] {1} (E);
\draw[->] (E) to node[left] {1} (A);

%bottom right
\node (A2) at (7,0) {$2^+$};
\node (B2) at (10,0) {$2^-$};
\node (C2) at (10,2) {$1^-$};
\node (D2) at (8.5,3.5) {$0$};
\node (E2) at (7,2) {$1^+$};
\draw[->] (A2) to node[below] {1} (B2);
\draw[->] (B2) to node[right] {0.5} (C2);
\draw[->] (C2) to node[above] {0.5} (A2);
\draw[->] (C2) to node[above] {\,\,\,\,0.5} (D2);
\draw[->] (D2) to node[above] {1} (E2);
\draw[->] (E2) to node[left] {0.5} (A2);
\draw[->] (E2) to node[above] {0.5} (C2);
\draw[->] (B2) to[bend right=80] node[above] {\,\,\,\,\,\,0.5} (D2);

%top left
\node (A3) at (-1,7) {$(2^+,1, 0, ..., 0)$};
\node (B3) at (3,7) {$(2^-,1, 0, ..., 0)$};
\node (C3) at (3,9) {$(1^-,0, ..., 0)$};
\node (D3) at (3,11) {$(2^+,0, ..., 0)$};
\node (E3) at (1,12.5) {$(2^-,0, ..., 0)$};
\node (F3) at (-1,11) {$(0, ..., 0)$};
\node (G3) at (-1,9) {$(1^+,0, ..., 0)$};
\draw[->] (A3) to node[below] {0.5} (B3);
\draw[->] (B3) to node[right] {0.5} (C3);
\draw[->] (C3) to node[right] {0.5} (D3);
\draw[->] (C3) to node[above] {0.5} (F3);
\draw[->] (D3) to node[above] {\,\,\,\,0.5} (E3);
\draw[->] (E3) to node[above] {0.5\,\,\,\,} (F3);
\draw[->] (F3) to node[left] {1} (G3);
\draw[->] (G3) to node[below] {0.5} (C3);
\draw[->] (G3) to node[left] {0.5} (A3);

%top right
\node (A2) at (7,7) {$2^+$};
\node (B2) at (10,7) {$2^-$};
\node (C2) at (10,9) {$1^-$};
\node (D2) at (8.5,10.5) {$0$};
\node (E2) at (7,9) {$1^+$};
\draw[->] (A2) to node[below] {1} (B2);
\draw[->] (B2) to node[right] {0.5} (C2);
\draw[->] (C2) to node[above] {0.5} (A2);
\draw[->] (C2) to node[above] {\,\,\,\,0.5} (D2);
\draw[->] (D2) to node[above] {1} (E2);
\draw[->] (E2) to node[left] {0.5} (A2);
\draw[->] (E2) to node[above] {0.5} (C2);
\draw[->] (B2) to[bend right=80] node[above] {\,\,\,\,\,\,0.5} (D2);

\node at (5,13) {$Q = 6$};
\node at (5,4) {$Q = 5$};
\node at (1,14) {event-based};
\node at (8.5,14) {location-based};
\end{tikzpicture}}
		\caption{Example for step two in the proof of Theorem \ref{T1}}
		\label{fig:ProofT1}
	\end{figure}
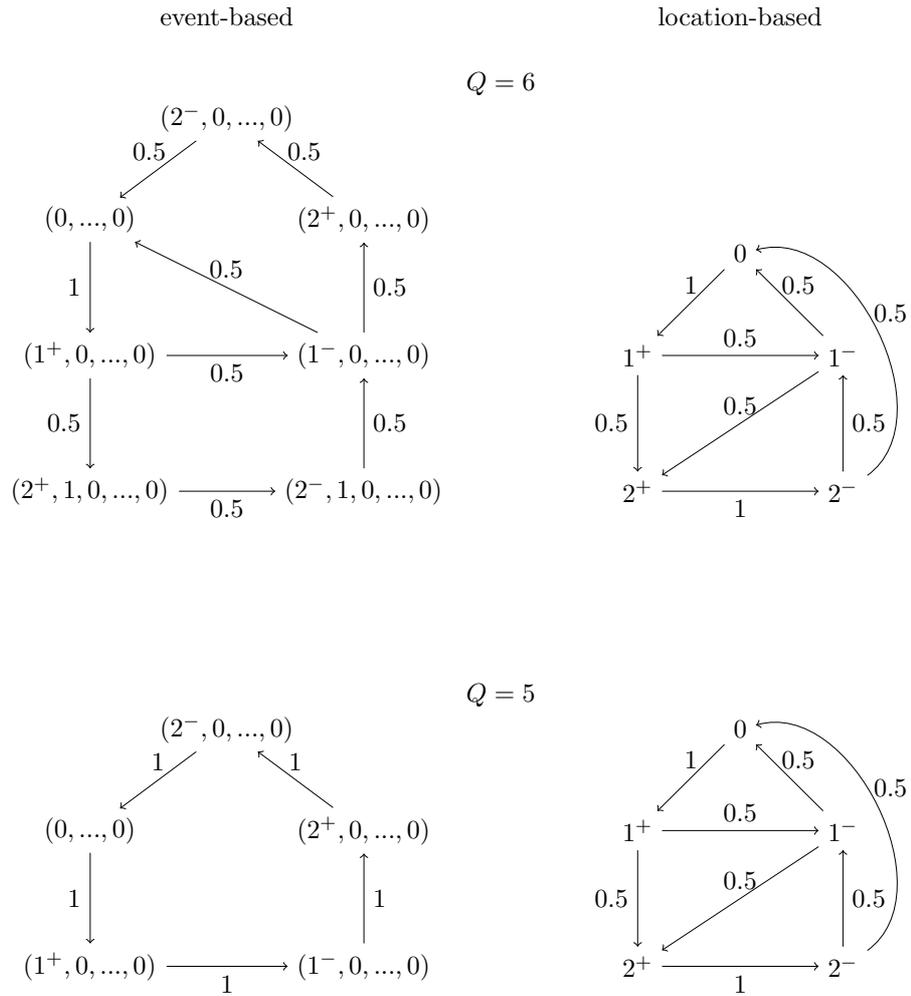
	Figure \ref{fig:ProofT1} shows an instance with two requests, $q_1 = q_2 = 3$, and one vehicle. Furthermore, we assume that time windows are not binding. The figure shows feasible LP solutions for the EB formulation on the left side and for the LB formulation on the right side. Note that only arcs with positive value for $x_{(v,w)}$ and $\bar{x}_{ij}$, respectively, are included in the figure. In the upper part of the figure, the vehicle capacity is $6$, i.e., both customers can be transported simultaneously. The graphs of both approaches show the solutions
	\[
	\begin{aligned}
		&0 \rightarrow 1^+ \rightarrow 2^+ \rightarrow 2^- \rightarrow 1^- \rightarrow 0 \\
		&0 \rightarrow 1^+ \rightarrow 1^- \rightarrow 2^+ \rightarrow 2^- \rightarrow 0 \\
	\end{aligned}
	\]
	both with a weight of 0.5. If we reduce the vehicle capacity to 5, the first solution is not integer feasible any more. While the event-based graph pictures this fact (lower left part of the figure), the solution in the LB formulation does not change. The reason is that capacity constraints are still fulfilled if the first solution is weighted with 0.5 and $Q_0 = 0, Q_{1^+} = 3$, $Q_{1^-} = 0$, $Q_{2^+} = 5$, and $Q_{2^-} = 2$ as the following evaluation shows: 
	\[
	\begin{aligned}
		& 										&&Q_{i}+ \mathcal{I}_{\{ j \in P \}}q_j - Q(1-x_{ij}) 		&&&\leq Q_{j} \\
		&0 \rightarrow 1^+  	&&0 + 3 - 5 \cdot (1-0.5-0.5) = 3 			&&& \leq 3 \\
		&1^+ \rightarrow 1^- 	&&3 - 3 - 5 \cdot 0.5 = -2.5 	&&&\leq 0 \\
		&1^+ \rightarrow 2^+ 	&&3 + 3 - 5 \cdot 0.5 = 3.5 	&&&\leq 5 \\
		&1^- \rightarrow 0  	&&0 + 0 - 5 \cdot 0.5 = -2.5 	&&& \leq 0 \\
		&1^- \rightarrow 2^+ 	&&0 + 3 - 5 \cdot 0.5 = 0.5 	&&&\leq 5 \\
		&2^+ \rightarrow 2^- 	&&5 - 3 - 5 \cdot (1-0.5-0.5) = 2 			&&&\leq 2 \\
		&2^- \rightarrow 0  	&&2 + 0 - 5 \cdot 0.5 = -0.5 	&&&\leq 0 \\
		&2^- \rightarrow 1^- 	&&2 - 3 - 5 \cdot 0.5 = -3.5 	&&&\leq 0 \\
	\end{aligned}
	\]
	In all other cases, $\bar{x}_{ij} = 0$. Since $Q_j \geq \max \{ 0,\mathcal{I}_{\{ j \in P \}}q_j \}$, it holds that $Q_i -Q \leq 0 \leq  Q_j - \mathcal{I}_{\{j\in P\}}q_j$ and the inequality is fulfilled in these cases. Thus, there is an LP solution in the LB formulation which is not feasible for the EB formulation.
	\item It remains to show that these LP solutions can be optimal. 
	\begin{table}[htbp]
		\begin{center}
			\begin{tabular}{cccccc}
				&	0		&	$1^+$	& $2^+$	&	$1^-$	&	$2^-$ \\
				0			& 0 	&	1			&	1			&	10		&	10		\\
				$1^+$	&	1		& 0			& 2			&	9			& 9			\\
				$2^+$	&	1		& 2			& 0			&	9			& 9			\\
				$1^-$	&	10	& 9			& 9			&	0			&	1			\\
				$2^-$	&	10	& 9			& 9			&	1			&	0			\\
			\end{tabular}
			\caption{Example for step three in the proof of Theorem \ref{T1}}
			\label{tab:T1}
		\end{center}
	\end{table}
	With the travel costs in Table \ref{tab:T1} the solution $0 \rightarrow 1^+ \rightarrow 2^+ \rightarrow 2^- \rightarrow 1^- \rightarrow 0$ has an objective value of 23 and the solution $0 \rightarrow 1^+ \rightarrow 1^- \rightarrow 2^+ \rightarrow 2^- \rightarrow 0$ an objective value of 38. Moreover, swapping $1^+$ and $2^+$ or $2^-$ and $1^-$ in the first solution or swapping the indices of customers 1 and 2 leads to the same objective values due to symmetry. Therefore, the LB model will use as much as possible of the first solution and its symmetric ones (such that capacity constraints stay fulfilled) while the EB model uses the second solution. Thus, $Z_{LAEB}^{rel} > Z_{LB}^{rel}$ holds.
\end{enumerate}
%\end{Proof}

\section{Implementation of Infeasible Path Constraints }\label{ImplIP}
For the implementation of inequalities \eqref{IP1} and \eqref{IP2} we first compute $B_{vw}^{LB}$ and $B_{vw}^{UB}$ for all $(v,w) \in A$. To implement inequalities \eqref{IP1}, for all events $v$, we loop over all successors $u$ (in decreasing order of $B_{vu}^{UB}$), determine the first predecessor $w$ of $v$ with $B_{vu}^{UB}<B_{wv}^{LB}$, and add all valid inequalities of type \eqref{IP1}, i.e., for each $v$ we add the inequality with variables $x_{(v,u')}$ such that $B_{vu'}^{UB} \leq B_{vu}^{UB}$ and variables $x_{(w',v)}$ such that $B_{w'v}^{LB} \geq B_{wv}^{LB}$. 
%It is useful to start with the successors, as there are up to $|R|+Q$ of them but up to $2|R|$ predecessors. 
Note that there is a dominance if for two successors $\bar{u}^1$ and $\bar{u}^2$ of the decreasing ordered list with $B_{v\bar{u}^1}^{UB} \geq B_{v\bar{u}^2}^{UB}$  the set of arcs $(w,v)$ with $B_{v\bar{u}^1}^{UB}<B_{wv}^{LB}$ and $B_{v\bar{u}^2}^{UB}<B_{wv}^{LB}$, respectively, are identical. Then, the valid inequality \eqref{IP1} for $\bar{u}^2$ is equal or weaker to the one for $\bar{u}^1$ and therefore not added.

For the implementation of \eqref{IP2} we consider the case $\bar{m}=4$. For all arcs $(u^2,u^3) \in A\setminus \left(\delta^\text{in}\left((0, \ldots,0)\right) \cup \delta^\text{out}\left((0, \ldots,0)\right)\right)$, we loop over all successors $u^4$ of $u^3$ (in decreasing order of $B_{u^3u^4}^{UB}$) and determine the first predecessor $u^1$  (in increasing order of $B_{u^1u^2}^{LB}$) of $u^2$ for which the path $u^1 \rightarrow u^2 \rightarrow u^3 \rightarrow u^4$ is infeasible. We can apply the same dominance criterion like in the implementation of \eqref{IP1} here. Again, we add all variables $x_{(u^3,w)}$ with $B_{u^3w}^{UB} \leq B_{u^3u^4}^{UB}$, $w \neq u^4$, and all variables $x_{(v,u^2)}$ with $B_{vu^2}^{LB} \geq  B_{u^1u^2}^{LB}$, $v \neq u^1$, to the left side of the inequality. We loop over all requests $i\in R$ and check whether events $u^2_{i^+}$ and $u^3_{i^-}$ and arcs $(u^2,u^2_{i^+})$ and $(u^3_{i^-},u^3)$ exist. If so, we add $\frac{1}{2} \cdot x_{(u^2,u^2_{i^+})}$ and $\frac{1}{2} \cdot x_{(u^3_{i^-},u^3)}$ to the inequality.

\section{Multi-thread Computational Results}
The following table shows the numerical results on the a-X and b-X instances solved with the LAEB model and the proposed preprocessing steps and without limiting the execution of CPLEX to one thread. In this case, 12 threads were used.
\setlength{\tabcolsep}{3pt}
\sisetup{table-align-text-post = false, table-number-alignment=center,table-format = 4.1}
\begin{table}[ht]
	\centering\footnotesize
	\begin{threeparttable}		
		\begin{tabular}{cSS[table-format = 1.1]Sc@{\hspace{2mm}}cSS}
			\toprule
			$\text{Inst.}$ &  $\text{Obj.}$  & $\text{Gap}$  & $\text{CPU}$  & 	& $\text{Inst.}$ &  $\text{Obj.}$ &  $\text{CPU}$  \\
			\midrule
			a2-16-X &278,2&&0,13&&b2-16-X&282,5&0,12 \\
			a2-20-X &330,7&&0,03&&b2-20-X&323,6&0,06 \\
			a2-24-X &389,1&&0,09&&b2-24-X&412,3&0,03 \\
			a3-18-X &272,7&&0,09&&b3-18-X&290,4&0,08 \\
			a3-24-X &289,6&&0,34&&b3-24-X&363,7&0,12 \\
			a3-30-X &452,8&&0,14&&b3-30-X&504,3&0,14 \\
			a3-36-X &501,0&&0,20&&b3-36-X&565,9&0,07 \\
			a4-16-X &235,2&&0,20&&b4-16-X&289,9&0,03 \\
			a4-24-X &359,4&&0,15&&b4-24-X&347,0&0,27 \\
			a4-32-X &447,3&&0,76&&b4-32-X&491,0&0,10 \\
			a4-40-X &509,0&&1,75&&b4-40-X&628,3&0,21 \\
			a4-48-X &620,3&&15,88&&b4-48-X&627,4&1,18 \\
			a5-40-X &464,0&&2,17&&b5-40-X&585,1&1,39 \\
			a5-50-X &621,9&&15,94&&b5-50-X&708,8&0,48 \\
			a5-60-X &745,4&&16,68&&b5-60-X&851,9&1,72 \\
			a6-48-X &572,5&&101,85&&b6-48-X&691,6&0,54 \\
			a6-60-X &757,9&&100,03&&b6-60-X&841,6&1,04 \\
			a6-72-X &868,4&&4368,93&&b6-72-X&930,3&6,08 \\
			a7-56-X &663,5&&102,02&&b7-56-X&787,9&6,50 \\
			a7-70-X &815,3&&100,12&&b7-70-X&865,3&2,77 \\
			a7-84-X &950,6&&6396,27&&b7-84-X&1141,2&50,73 \\
			a8-64-X &701,2&&191,76&&b8-64-X&818,3&3,04 \\
			a8-80-X &880,3&&4579,24&&b8-80-X&998,3&2,73 \\
			a8-96-X &1118,5&2.3 \%&2h&&b8-96-X&1137,7&147,31 \\
			\cmidrule(lr){1-1} \cmidrule(lr){4-4} \cmidrule(lr){8-8}
			Total &              &&  \text{23196} & & & & \text{226} \\
			\bottomrule		
		\end{tabular}
		\caption{Results for the a-X and b-X benchmark instances on 12 threads.}
		\label{tab:12threads_aXbX}
	\end{threeparttable}
\end{table}

d em
%%%%%%%%%%%%%%%%%
\end{document}